\documentclass[11pt]{article}

\usepackage{amssymb,amsmath,latexsym}
\usepackage[left=1in, right=1in,top=1in, bottom=1in]{geometry}
\usepackage{amsfonts}
\usepackage{mathrsfs}
\usepackage{cases}
\usepackage{bm}
\usepackage{indentfirst}
\usepackage{color}
\usepackage{ifpdf}
\usepackage{graphicx}
\usepackage{psfrag}
\usepackage{dsfont}
\usepackage{thmtools}
\usepackage{thm-restate}
\usepackage{cleveref}
\usepackage{enumerate}

\usepackage[
pdfauthor={CS},
pdftitle={Halin Graphs},
pdfstartview=XYZ,
bookmarks=true,
colorlinks=true,
linkcolor=blue,
urlcolor=blue,
citecolor=blue,
bookmarks=true,
linktocpage=true,
hyperindex=true
]{hyperref}

\newcommand{\xiaowuhao}{\fontsize{9pt}{\baselineskip}\selectfont}

\newtheorem{THM}{\textbf{Theorem}}[section]
\newtheorem{DEF}{\textbf{Definition}}
\newtheorem{LEM}{\textbf{Lemma}}[section]

\newtheorem{PRO}{\textbf{Proposition}}

\newtheorem{CLA}{\textbf{Claim}}[subsection]

\newcommand{\qed}{\hfill $\square$\vspace{1mm}}

\newcommand{\ve}{\varepsilon }
\newcommand{\nb}{\noindent \textbf }

\newcommand{\pf}{\textbf{Proof}.\quad}

\begin{document}
\title{Dirac's Condition for Spanning Halin Subgraphs}
\author{Guantao Chen$^{\dag}$ and Songling Shan$^{\ddag}$\\
{\xiaowuhao $^{\dag}$  Georgia State University, Atlanta, GA\,30303 }\\
{\xiaowuhao  $^{\ddag}$ Vanderbilt University, Nashville, TN\,37240}}

\date{\today}
\maketitle
\abstract{Let $G$ be an $n$-vertex graph with $n\ge 3$.
A classic result of Dirac from 1952 asserts that $G$ is hamiltonian
if $\delta(G)\ge n/2$. Dirac's theorem is one of the most
influential results in the study of hamiltonicity and by now
there are many related known results\,(see, e.g., J. A. Bondy, Basic Graph Theory:
Paths and Circuits, Chapter 1 in: {\it Handbook of Combinatorics Vol.1}).
A {\it Halin graph} is a planar graph consisting of two edge-disjoint subgraphs:
a spanning tree of at least 4 vertices and with no vertex of degree 2,
and a cycle induced on the set of the leaves of the spanning tree.
Halin graphs possess rich hamiltonicity properties
such as being hamiltonian, hamiltonian connected, and almost pancyclic.
As a continuous ``generalization'' of Dirac's theorem,
in this paper, we show that
there exists a positive integer $n_0$ such that
any graph $G$ with $n\ge n_0$ vertices and
$\delta(G)\ge (n+1)/2$
contains a spanning Halin subgraph.
In particular, it contains a spanning Halin subgraph which is also
pancyclic.
}
\medskip

\section{Introduction}
A classic theorem of Dirac~\cite{Dirac-thm} from 1952 asserts that every graph on $n$ vertices
with minimum degree at least $n/2$ is hamiltonian if $n\ge 3$. Following Dirac's result, numerous
results on hamiltonicity properties on graphs with restricted degree conditions have been obtained\,
(see, for instance, \cite{Ron-gold-hamiltonian-survey, MR1373655}). Traditionally,  under similar conditions,  results for a graph being hamiltonian, hamiltonian-connected, and pancyclic are obtained separately. We may ask, under certain conditions,  if it is possible to uniformly show a graph possessing several hamiltonicity properties. The work on finding the square of a hamiltonian cycle
in a graph can be seen as an attempt in this direction. However, it requires minimum  degree of $2n/3$
for an $n$-vertex
graph $G$ to contain the square of a hamiltonian cycle,
 for examples, see~\cite{Degreesum-H2, Posa-H2conjecture,MR1399673, Komlos,  1412.3498}.
Although the minimum degree condition of $2n/3$ for having the square of a hamiltonian cycle
is almost optimal for the embedding result it implies:
 Aigner-Brandt Theorem~\cite{MR1223891} that any $n$ vertex graph with minimum degree at least $(2n-1)/3$
contains  all possible
 graphs of order at most $n$ and maximum degree at most 2, it is a ``waste'' for
 using   the square of a hamiltonian cycle in obtaining hamiltonicity properties.
 For bipartite graphs, finding the existence of a spanning ladder is a way of simultaneously showing the graph having many hamiltonicity properties~(see~\cite{2-factor-bipartite, MR2646098}).    In this paper,  we introduce another approach of uniformly showing the possession of several hamiltonicity properties in a graph: we show the existence of a spanning {\it Halin graph} in a graph under a
 given minimum degree condition.


A tree with no vertex of degree 2 is called a {\it homeomorphically irreducible tree}\,(HIT).
A {\it Halin graph }  $H=T\cup C$  is a simple  planar graph
 consisting of  a  HIT $T$  with at least 4 vertices  and a
 cycle $C$  induced by   the set of leaves of $T$.
The HIT $T$ is called the {\it underlying tree} of $H$.
A wheel graph  is an example of a Halin graph, where
the underlying tree is a star.
Halin constructed Halin graphs in~\cite{Halin-halin-graph} for the study of minimally 3-connected graphs.
Lov\'asz and Plummer named such graphs as Halin graphs
in their study of planar bicritical graphs~\cite{LP-Halingraph-Con},
which are planar graphs having a 1-factor after deleting any two vertices.
Intensive researches  have been done on Halin graphs.
Bondy~\cite{Bondy-pancyclic} in 1975 showed that a Halin graph is hamiltonian.
In the same year, Lov\'asz and Plummer~\cite{LP-Halingraph-Con} showed that not only a Halin graph
itself is hamiltonian, but each of the subgraph obtained  by
deleting a vertex is hamiltonian.
In 1987, Barefoot~\cite{Barefoot} proved that Halin graphs are hamiltonian-connected,
i.e., there is a hamiltonian path connecting any two vertices of the graph.
Furthermore,  it was proved that each edge of a Halin graph is contained
in a hamiltonian cycle and is avoided by another~\cite{Skupie-uniformly-hamiltonian}.
Bondy and Lov\'asz~\cite{Almost-pancyclic-halin}, and Skowro\'nska~\cite{pancyclicity-Halin-graphs}, independently, in 1985,
showed that a Halin graph is almost pancyclic and is pancyclic if the underlying tree has
no vertex of degree 3, where an $n$-vertex graph is {\it almost pancyclic} if it contains cycles of
length from 3 to $n$ with the possible exception
  of a single even length, and  is {\it pancyclic} if it contains cycles of
length from 3 to $n$.  Some problems that are NP-complete for general graphs have been shown to
be polynomial time solvable for Halin graphs.
For example,
Cornu\'ejols, Naddef, and  Pulleyblank\,\cite{Cornuejols-halin} showed that
in a Halin graph, a hamiltonian cycle can be found in polynomial time.
It seems so promising to show the existence of a spanning Halin subgraph in
a given graph in order to show that the graph possesses many hamiltonicity properties.
But, nothing comes for free,
it is NP-complete to determine whether a graph contains a (spanning) Halin graph~\cite{Halin-NP}.

Despite all these nice properties of Halin graphs mentioned above, the problem of determining whether
a graph contains a spanning Halin subgraph has not yet well studied except a conjecture
proposed by Lov\'asz and Plummer~\cite{LP-Halingraph-Con} in 1975.
The conjecture states that {\it every 4-connected plane
triangulation contains a spanning Halin subgraph}\,(disproved recently~\cite{Disprove-LP-Con}).
In this paper, we investigate the minimum degree condition for implying the existence of
a spanning Halin subgraph in a graph, and thereby giving another approach  for
uniformly showing the possession of several hamiltonicity properties in a graph under a
given minimum degree condition.  We obtain the following result.

\begin{THM}\label{main result}
There exists $n_0>0$ such that for  any graph $G$ with $n\ge n_0$
vertices,   if $\delta(G)\ge (n+1)/2$, then $G$ contains a spanning Halin subgraph.
In particular, it contains a spanning Halin subgraph which is also
pancyclic.
\end{THM}

Note that an $n$-vertex graph with minimum degree at least $(n+1)/2$
is 3-connected if $n\ge 4$. Hence, the minimum degree condition in Theorem~\ref{main result}
implies the 3-connectedness, which is a necessary condition for a graph to contain a
spanning Halin subgraph, since every Halin graph is 3-connected.
A Halin graph contains a triangle, and bipartite graphs are triangle-free.
Hence, $K_{\lfloor \frac{n}{2}\rfloor,\lceil \frac{n}{2}\rceil}$  contains
no spanning Halin subgraph.  For $n$ even, the graph obtained
from two copies of $K_{\frac{n}{2}+1}$ by gluing them together
on an edge is 2-connected, so it has no  spanning Halin subgraph.
Both these two graphs have minimum degree at most $n/2$.
We see that the minimum degree condition
in Theorem~\ref{main result} is best possible.

Theorem~\ref{main result} is proved for large graphs. It
might be very hard for obtaining a same result for all graphs,
as when constructing a Halin graph in general,
we may need to find its underlying tree first. The minimum degree
condition suffices for the existence of  a such tree $T$ in $G$
(in fact, it was showed that an $n$-vertex graph with minimum degree at least
$4\sqrt{2n}$ contains a spanning tree with no vertex of degree 2~\cite{MR1053607}).
However, the hardness lies in finding a cycle $C$ spanning on
the set of the leaves of $T$ so that $T\cup C$ is planar.
In other words, when $T$ is fixed, we have to find a cycle $C$
in $G$ passing through a set of given vertices in some particular  order.
The other way of finding a spanning Halin graph $H$ is
to find a spanning  subgraph which contains $H$.
For example, spanning structures close to
ladder structures\,(e.g., graphs $H_1$ to $H_5$ as defined in next section).
Particularly, the square of a hamiltonian cycle contains
$H_1$ or $H_2$ as a spanning subgraph, so it contains a spanning Halin subgraph.
But the disadvantage of using ``uniform'' structures as $H_i$
is that it makes it hard for  constructing them
``manually''.  Nevertheless, we still suspect that
$(n+1)/2$ is the right condition for all graphs
to contain a spanning Halin subgraph.


\section{Notation and definitions}
We consider simple and finite graphs only.
Let $G$ be a graph. Denote by $V(G)$ and
$E(G)$ the vertex set and edge set of $G$, respectively,
and by $e(G)$ the cardinality of $E(G)$. We denote by
$\delta(G)$ the minimum degree of $G$ and by $\Delta(G)$
the maximum degree. Let $v\in V(G)$ be a vertex and $S\subseteq V(G)$
a subset. Then $G[S]$ is the subgraph of $G$ induced by $S$.
Similarly, $G[F]$ is the subgraph induced by $F$ if
$F\subseteq E(G)$.
The notation $\Gamma_G(v,S)$ denotes the set of neighbors of
$v$ in $S$, and $deg_G(v,S)=|\Gamma_G(v,S)|$.
We let $\Gamma_{\overline{G}}(v,S)=S-\Gamma_G(v,S)$ and
$deg_{\overline{G}}(v,S)=|\Gamma_{\overline{G}}(v,S)|$.
Given another set $U\subseteq V(G)$,
define $\Gamma_G(U,S)=\bigcap_{u\in U}\Gamma_G(u,S)$, $deg_G(U,S)=|\Gamma_G(U,S)|$, and $N_G(U,S)=\bigcup_{u\in U}\Gamma_G(u,S)$.
When $U=\{u_1,u_2,\cdots, u_k\}$,  we may  write $\Gamma_G(U,S)$,
$deg_G(U,S)$, and $N_G(U,S)$ as $\Gamma_G(u_1,u_2,\cdots, u_k,S)$,  $deg_G(u_1,u_2,\cdots, u_k,S)$, and $N_G(u_1,u_2,\cdots, u_k,S)$,
respectively, in specifying the vertices in $U$.  When $S=V(G)$, we only write $\Gamma_G(U)$, $deg_G(U)$, and $N_G(U)$.
Let $U_1,U_2 \subseteq V(G)$ be two disjoint subsets.
Then  $\delta_G(U_1,U_2)=\min\{deg_G(u_1,U_2)\,|\, u_1\in U_1\}$ and
$\Delta_G(U_1,U_2)=\max\{deg_G(u_1,U_2)\,|\, u_1\in U_1\}$.
Notice that the notation $\delta_G(U_1,U_2)$ and $\Delta_G(U_1,U_2)$
are not symmetric with respect to $U_1$ and $U_2$.
We denote by $E_G(U_1,U_2)$ the set of edges with one end in $U_1$ and the other in $U_2$,
the cardinality of $E_G(U_1,U_2)$ is denoted by  $e_G(U_1,U_2)$.
We may omit the index $G$ if there is no risk of confusion.
Let $u,v\in V(G)$
be two vertices. We write $u\sim v$ if $u$ and $v$  are adjacent.
A path connecting  $u$ and $v$  is called
a  $(u,v)$-path.
If $G$ is a bipartite graph with partite sets
$A$ and $B$, we denote $G$ by $G(A,B)$ in emphasizing the two partite sets.

In constructing Halin graphs, we use ladder graphs and a class of ``ladder-like''
graphs as substructures. We give the description of these graphs below.

\begin{DEF}\label{ladder}
An $n$-ladder $L_n=L_n(A,B)$ is a balanced bipartite graph with
$A=\{a_1,a_2,\cdots, a_n\}$ and $B=\{b_1,b_2,\cdots, b_n\}$
such that $a_i\sim b_j$ iff $|i-j|\le 1$. We call $a_ib_i$ the $i$-th
rung of $L_n$. If $2n (mod \,\,4) \equiv 0$, then we call each of the shortest
$(a_1,b_n)$-path $a_1b_2a_3b_4\cdots a_{n-1}b_n$ and $(b_1,a_n)$-path $b_1a_2b_3a_4\cdots b_{n-1}a_n$ a side of $L_n$;
and if  $2n (mod \,\,4) \equiv 2$, then  we call each of the shortest
$(a_1,a_n)$-path $a_1b_2a_3b_4\cdots a_{n-1}b_{n-1}a_n$ and $(b_1,b_n)$-path $b_1a_2b_3a_4\cdots b_{n-2}a_{n-1}b_n$ a side of $L_n$.
\end{DEF}

%

Let $L$ be a ladder with $xy$ as one of its rungs.   For an edge $gh$,
we say $xy$ and $gh$ are {\it adjacent} if  $x\sim g, y\sim h$
or $x\sim h, y\sim g$.  Suppose
$L$ has its first rung as $ab$ and its last rung as $cd$,
we denote $L$ by $ab-L-cd$ in specifying the two rungs,
and we always assume that the distance between $a$ and $c$ and thus between
$b$ and $d$
is $|V(L)|/2-1$\,(we make this assumption for being convenient in constructing
other graphs based on ladders).  Under this assumption, we
denote $L$ as $\overrightarrow{ab}-L-\overrightarrow{cd}$.
Let $A$ and $B$ be two disjoint vertex sets.
We say the rung $xy$ of $L$ is {\it contained} in $A\times B$
if either $x\in A, y\in B$ or $x\in B, y\in A$.  Let $L'$ be another
ladder vertex-disjoint with $L$. If the last rung of $L$ is adjacent to
the first rung of $L'$, we write $LL'$ for the new ladder
obtained by concatenating $L$ and $L'$. In particular, if $L'=gh$
is an edge, we write $LL'$ as $Lgh$.

\begin{figure}[!htb]
\psfrag{L_4}{$L_4$} \psfrag{H_1}{$H_1$}
\psfrag{H_2}{$H_2$}
\psfrag{H_3}{$H_3$}
\psfrag{H_4}{$H_4$}
\psfrag{H_5}{$H_5$}
\psfrag{a1}{$a_1$} \psfrag{a2}{$a_2$} \psfrag{a3}{$a_3$} \psfrag{a4}{$a_4$}
\psfrag{b1}{$b_1$} \psfrag{b2}{$b_2$} \psfrag{b3}{$b_3$} \psfrag{b4}{$b_4$}
\psfrag{x}{$x$} \psfrag{y}{$y$} \psfrag{z}{$z$} \psfrag{u}{$u$} \psfrag{w}{$w$}
\psfrag{L_4}{$L_4$} \psfrag{T_1}{$T_1$}
\psfrag{T_2}{$T_2$}
\psfrag{T_3}{$T_3$}
\psfrag{T_4}{$T_4$}
\psfrag{T_5}{$T_5$}
\begin{center}
  \includegraphics[scale=0.38]{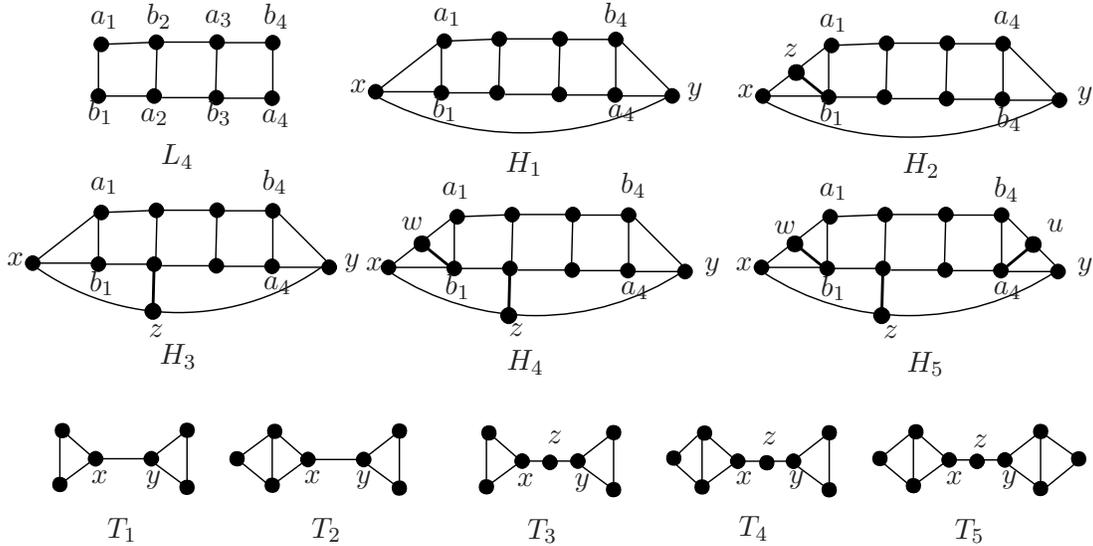}\\
\end{center}
\vspace{-3mm}
  \caption{{\small $L_4$, $H_i$  constructed from $L_4$, and $T_i$
  associated with $H_i$ for each $i=1,2,\cdots, 5$
}}\label{cydder}
\end{figure}

We now define  five types of ``ladder-like''  graphs, call them
$H_1, H_2, H_3, H_4$ and $H_5$, respectively.
Let  $L_{n}$ be a ladder
with $a_1b_1$  and $a_{n}b_{n}$  as the first and
last rung, respectively, and
$x,y,z,w, u$ be five new vertices. Then each of $H_i$
is obtained from $L_{n}$ by adding some
specified vertices and edges as follows. Additionally, for each $i$ with $1\le i\le 5$,
we define a graph $T_i$ associated with $H_i$. A depiction of
a ladder $L_4$, $H_1,H_2,H_3,H_4,H_5$
constructed from $L_4$, and the graph $T_i$ associated with $H_i$
is given in Figure~\ref{cydder}.

\begin{itemize}
  \item [$H_1$: ] Adding two new vertices $x, y$ and the edges $xa_1,xb_1, ya_{n}, yb_{n}$
and $xy$.

Let $T_1=H_1[\{x,y,a_1,b_1, a_n, b_n\}]$.
  \item [$H_2$: ] Adding three new vertices $x,y,z$ and the edges $za_1,zb_1, xz,xb_1,  ya_{n}, yb_{n}$
and $xy$.

Let $T_2=H_2[\{x,y,z, a_1,b_1, a_n, b_n\}]$.
  \item [$H_3$: ] Adding three new vertices $x,y,z$ and the edges $xa_1,xb_1,  ya_{n}, yb_{n}, xz, yz$,  and
  either $za_{i}$ or $zb_i$ for
some $1\le i \le n$. Note that $H_2$ is a special case of $H_3$ with $i=1$ or $n$.

 Let $T_3=H_3[\{x,y,z,a_1,b_1, a_n, b_n\}]$.
  \item [$H_4$: ] Adding four  new vertices $x,y,z, w$ and the edges $wa_1,wb_1,  xw, xb_1, ya_{n}, yb_{n},xz,yz$,
  and either $za_{i}$ or $zb_i$ for
some $1\le i \le n$ such that $a_i$ or $b_i$ is a vertex on the side of $L$ which has $b_1$ as one end.

Let $T_4=H_4[\{x,y,z,w,a_1,b_1, a_n, b_n\}]$.
 \item [$H_5$: ] Adding five new vertices $x,y,z, w,u$.

 If  $2n (mod\,\, 4) \equiv 2$, adding the edges $wa_1,wb_1,  xw, xb_1, ua_{n}, ub_{n}, yu, yb_{n}, xz,yz$,  and
 either $ za_{i}$ or $zb_i$ for
some $1\le i \le n$ such that $a_i$ or $b_i$ is a vertex on the shortest $(b_1,b_n)$-path in $L$;

and if $2n (mod \,\,4) \equiv 0$, adding
the edges $wa_1,wb_1,  xw, xb_1, ua_{n}, ub_{n},  yu, ya_{n}, xz, yz$, and either $ za_{i}$ or $zb_i$ for
some $1\le i \le n$ such that $a_i$ or $b_i$ is a vertex on the shortest $(b_1,a_n)$-path in $L$.

The graph obtained from $H_5$ by deleting the vertex $z$ and adding the edge $xy$ is identical with
$H_4$ with $i=n$.

Let $T_5=H_5[\{x,y,z,w,u,a_1,b_1, a_n, b_n\}]$.

\end{itemize}

Let $i=1,2,\cdots, 5$.   Notice that each of $H_i$ is a Halin graph,
and the graph obtained from $H_5$ by deleting the vertex $z$ and adding the edge $xy$
is also a Halin graph.
Except $H_1$, each $H_i$ has a unique underlying tree.
Notice also that $xy$ is an edge on the cycle
along the leaves of any  underlying tree  of $H_1$ or $H_2$.
For each $H_i$ and $T_i$,  call $x$
the {\it left end} and $y$ the {\it right end},
and call  a vertex of degree at least 3 in the underlying tree of
$H_i$
a {\it Halin constructible vertex}.
By analyzing the structure of $H_i$, we see that each internal vertex on a/the shortest $(x,y)$-path
in $H_i-xy$\,(for $i=1,2$) or $H_i-z$\,(for $i=3,4,5$)
is
a Halin constructible vertex.
Noting that any vertex in $V(H_1)-\{x,y\}$ can be a Halin constructible vertex.
We call $a_1b_1$ the head link of $T_i$ and $a_nb_n$ the tail link of $T_i$, and for each  of $T_3, T_4, T_5$, we call the vertex $z$
not contained in any triangles the {\it pendent vertex}. The notation of $H_i$ and $T_i$
are fixed hereafter.

Let $T\in \{T_1,\cdots, T_5\}$ be a subgraph of a graph $G$. Suppose that $T$ has head link $ab$, tail link  $cd$,
and possibly the pendent vertex $z$.
Suppose  $G-V(T)$ contains a spanning ladder $L$ with first rung
$c_1d_1$ and last rung $c_nd_n$ such that $c_1d_1$ is adjacent to $ab$,
$c_nd_n$ is adjacent to $cd$. Additionally,
if the pendent vertex $z$ of $T$ exists, then $z$ has
a neighbor $z'$,  which is an internal vertex on a shortest path between the two ends
of $T$ in the graph $abLcd\cup T-z$.
Then $abLcd\cup T$ or $abLcd\cup T \cup \{zz'\}$
is a  spanning Halin subgraph of $G$. This technique is  frequently used later on in constructing a
Halin graph. The following proposition  gives another way of constructing a Halin graph based on $H_1$ and $H_2$.

\begin{PRO}\label{Prop:Halin_H1_H2}
For $i=1,2$, let $G_i\in \{H_1,H_2\}$ with left end $x_i$ and right end $y_i$ be defined as above,
and let $u_i\in V(G_i)$  be a Halin constructible vertex, then $Q:=G_1\cup G_2-\{x_1y_1,x_2y_2\}\cup\{x_1x_2, y_1y_2,u_1u_2\}$
is  a Halin graph spanning on $V(G_1)\cup V(G_2)$. Let $\mathcal{Q}=\{Q\,|\, G_i\in \{H_1,H_2\} \}$
be the set of all graphs $Q$ constructed in this way.
Then any graph in  $\mathcal{Q}$ is pancyclic.
\end{PRO}

\pf For $i=1,2$, let $G_i$ be embedded in the plane, and let $T_{G_i}$ be a underlying plane tree of $G_i$.
Then $T':=T_{G_1}\cup T_{G_2}\cup \{u_1u_2\}$ is a homeomorphically irreducible tree
spanning on $V(G_1)\cup V(G_2)$. Moreover, we can draw the edge $u_1u_2$ such that
$T_{G_1}\cup T_{G_2}\cup \{u_1u_2\}$ is a plane graph. Since $G_i[E(G_i)-E(T_{G_i})-\{x_iy_i\}]$ is an
$(x_i,y_i)$-path spanning on the set of leaves of $T_{G_i}$ obtained by connecting the
leaves following the order determined by the embedding of $T_{G_i}$, we see that
$G_1[E(G_1)-E(T_{G_1})-\{x_1y_1\}]\cup G_2[E(G_2)-E(T_{G_2})-\{x_2y_2\}]\cup \{x_1x_2, y_1y_2\}$
is a cycle spanning on the set of leaves of $T'$ obtained by connecting the leaves
following the order determined by the embedding of $T'$. Thus
$Q$ is a Halin graph.

To see the pancyclicity of graphs in $\mathcal{Q}$, suppose that $H_1$
has $2n_1+2$ vertices and $H_2$ has $2n_2+3$ vertices.
It is easy to check that in $H_i$, there are $(x_i, y_i)$-paths
of length from $n_i+1$ to $|V(H_i)|-1$; in $H_1$, there are cycles of
length from 3 to $2n_1+2$; and in $H_2$, there are cycles of
length from 3 to $2n_2+3$.
Let $Q\in \mathcal{Q}$
such that $Q$ is constructed based on $H_1$ and $H_2$.
As $x_1x_2, y_1y_2\in E(Q)$,
then we see that $Q$ has all cycles of length from $n_1+n_2+4$ to $n_1+n_2+5=|V(Q)|$.
Together with cycles in $H_i$, we know that $Q$ contains
all cycles of length from 3 to $|V(Q)|$. The pancyclicity
of other graphs in  $\mathcal{Q}$ can be checked similarly.
\qed

%

%
%
%

\section{Proof of Theorem~\ref{main result}}

In this section, we prove Theorem~\ref{main result}. Following the standard setup of
proofs applying the Regularity Lemma, we divide the proof into  Non-extremal Case and
extremal cases.  For this purpose, we define the two extremal cases in the following.

Let $G$ be an $n$-vertex graph and $V$ its vertex set.
Let $0< \beta\le 1$ be a constant.   Let $W\subseteq V(G)$.
We say $W$ is an {\it approximate vertex-cut\/} of $G$ with parameter $\beta$
if there is a partition $V_1$ and $V_2$ of $V-W$
such that $e_G(V_1,V_2)\le \beta n^2$ and $\delta[G[V_i]]\ge \delta(G)-|W|-\beta n$
for each $i=1,2$.
The two extremal cases are defined as below.

{\nb{Extremal Case 1.}}  $G$ has an approximate vertex-cut of size at most  $5\beta n$ with parameter $\beta$.

{\nb{Extremal Case 2.}} There exists a partition $V_1\cup V_2$ of $V$
such that $|V_1|\ge (1/2-7\beta)n$ and $\Delta(G[V_1])\le \beta n$.

{\nb{Non-extremal Case.}} We say that an $n$-vertex graph with minimum degree at least $(n+1)/2$ is
in {\it Non-extremal Case\/ } if it is in neither    Extremal Case 1 nor  Extremal Case 2.

In Extremal Case 1,  we will show that $G$ contains a spanning Halin subgraph
isomorphic to a graph in $\mathcal{Q}$\,(defined in Proposition~\ref{Prop:Halin_H1_H2}).
In all other cases, we will construct a spanning subgraph of $G$ isomorphic to
$H_i$ for some $i\in \{1,2,3,4,5\}$. Note that each graph in $\mathcal{Q}$
and each
$H_i$ is pancyclic.  Hence, to prove Theorem~\ref{main result},
we only need to show the existence of the mentioned
graphs above.
The following three theorems deal with the Non-extremal Case and the
two extremal cases, respectively, and thus give a proof of
Theorem~\ref{main result}.

\begin{THM}\label{extremal1}
Suppose that $0<\beta \ll 1/(20\cdot17^3)$ and $n$ is a sufficiently large integer.
Let $G$ be a graph on $n$ vertices with $\delta(G)\ge (n+1)/2$.
If $G$ is in Extremal Case 1,  then $G$
contains a spanning Halin isomorphic to a graph in $\mathcal{Q}$\,(defined in Proposition~\ref{Prop:Halin_H1_H2})
as a subgraph.
\end{THM}


\begin{THM}\label{extremal2}
Suppose that $0<\beta \ll 1/(20\cdot17^3)$ and $n$ is a sufficiently large integer.
Let $G$ be a graph on $n$ vertices with $\delta(G)\ge (n+1)/2$.
If $G$ is in Extremal Case 2,  then $G$
contains a spanning Halin subgraph isomorphic to
some $H_i$, $i\in \{1,2,3,4,5\}$.
\end{THM}

\begin{THM}\label{non-extremal}
Let $n$ be a sufficiently large integer and $G$ an $n$-vertex graph with
$\delta(G)\ge (n+1)/2$. If $G$ is in the Non-extremal Case,  then $G$
has a spanning Halin subgraph isomorphic to $H_1$ or $H_2$.
\end{THM}

We need the following lemma  in each of the proofs of Theorems~\ref{extremal1} - \ref{extremal2}
in dealing with ``garbage'' vertices.

\begin{LEM}\label{absorbing}
Let $F$ be a graph  such that $V(F)$ is partitioned as  $S\cup R$.  Suppose that
(i) there are $|R|$ vertex-disjoint 3-stars\,(a 3-star is a copy of $K_{1,3}$) with the vertices in $R$
as their centers,
(ii) for any two vertices $u,v\in N(R,S)$, $deg(u,v,S)\ge 6|R|$,  and
(iii) for any three vertices $u,v, w\in N(N(R,S),S)$, $deg(u,v,w, S)\ge 7|R|$.
Then  there is a ladder spanning on $R$ and some other $7|R|-2$ vertices from $S$.
Particularly, the ladder has the vertices on its first and last rungs in $S$.
\end{LEM}

\begin{figure}[!htb]
\psfrag{w_1}{$w_1$} \psfrag{w_2}{$w_2$}
\psfrag{x_{11}}{$x_{11}$} \psfrag{x_{12}}{$x_{12}$} \psfrag{x_{8}}{$x_{8}$} \psfrag{x_{13}}{$x_{13}$}
\psfrag{x_{21}}{$x_{21}$} \psfrag{x_{22}}{$x_{22}$} \psfrag{x_{23}}{$x_{23}$}
\psfrag{y_{12}^1}{$y_{12}^1$} \psfrag{y_{23}^1}{$y_{23}^1$} \psfrag{y_{12}^2}{$y_{12}^2$}
\psfrag{y_{23}^2}{$y_{23}^2$}
\psfrag{u_1}{$u_1$}
\psfrag{z_1}{$z_1$}
\begin{center}
  \includegraphics[scale=0.5]{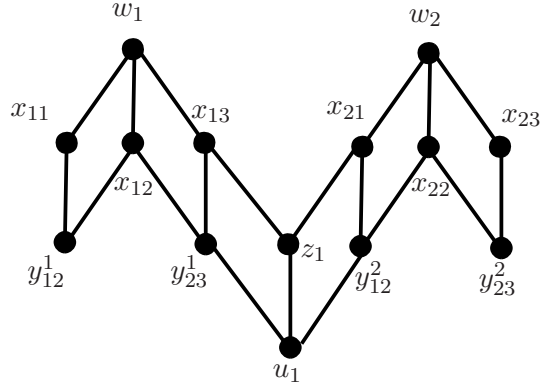}\\
\end{center}
\vspace{-3mm}
  \caption{Ladder $L$ of order 14}\label{insterL}
\end{figure}

\pf Let $R=\{w_1,w_2,\cdots, w_r\}$.
Consider first that $r=1$.  Choose $x_{11},x_{12},x_{13}\in \Gamma(w_1,S)$.
By (ii), there are distinct vertices $y^1_{12}\in \Gamma(x_{11},x_{12}, S)$
and $y^1_{23}\in \Gamma(x_{12},x_{13}, S)$.  Then the graph $L$ on $\{w_1,x_{11},x_{12},x_{13}, y^1_{12}, y^1_{23}\}$
with edges in
$$
\{w_1x_{11}, w_1x_{12}, w_1x_{13}, y^1_{12}x_{11}, y^1_{12}x_{12}, y^1_{23}x_{12}, y^1_{23}x_{13}\}
$$
is a ladder covering $R$ with $|V(L)|=6$. Suppose now $r\ge 2$.
By condition (i),
for each $i$ with $1\le i\le r$, there exist distinct  vertices $x_{i1}, x_{i2}, x_{i3}\in \Gamma(w_i, S)$.
By (ii), we choose distinct vertices $y_{12}^1, y_{23}^1,\cdots, y_{12}^r, y_{23}^r$
different from the existing vertices already chosen  such that $y_{12}^i\in \Gamma(x_{i1}, x_{i2}, S)$
and $y_{23}^i\in \Gamma(x_{i2}, x_{i3}, S)$ for each $i$, and  at the same  time, we chose
distinct vertices $z_1,z_2,\cdots, z_{r-1}$ from the unchosen  vertices in $S$
such that $z_i\in \Gamma(x_{i3}, x_{(i+1),1}, S)$ for each $1\le i\le r-1$.
Finally,  by (iii), choose distinct vertices $u_1, u_2,\cdots, u_{r-1}$
from the unchosen vertices in $S$ such that  $u_i\in \Gamma(y^i_{23}, y^{i+1}_{12}, z_i, S)$.
Let $L$ be the graph with
$$
V(L)=R\cup \{x_{i1}, x_{i2}, x_{i3},y_{12}^i, y_{23}^i, z_i, u_i, x_{r1}, x_{r2}, x_{r3}, y_{12}^r, y_{23}^r \,|\, 1\le i\le r-1 \}\quad \mbox{and}
$$
$E(L)$ consisting of the edges $w_rx_{r1}, w_rx_{r2},w_rx_{r3}, y_{12}^rx_{r1}, y_{12}^rx_{r2},y_{23}^rx_{r2},y_{23}^rx_{r3}$
and the edges indicated below for each $1\le i\le r-1$:
$$
w_i\sim x_{i1},x_{i2}, x_{i3};\, y_{12}^i\sim x_{i1}, x_{i2};\, y_{23}^i\sim x_{i2},x_{i3}; \,  z_i\sim x_{i3}, x_{i+1,1}; \, u_i\sim x_{i3}, x_{i+1,1}, z_i.
$$
It is easy to check that $L$ is a ladder covering $R$ with $|V(L)|=8r-2$. The ladder has its first and last rungs in $S$
is seen by its construction.
Figure~\ref{insterL}  gives a depiction  of  $L$ for $|R|=2$. \qed

We will also need the bipartite version of Lemma~\ref{absorbing}.
Since the proof is similar, we omit it.
\begin{LEM}\label{absorbing2}
Let $F(A,B)$ be a bipartite graph  such that $V(F)$ is partitioned as  $S\cup R$ with
$R\subseteq A$.  Suppose that
(i) there are $|R|$ vertex-disjoint 3-stars with the vertices in $R$
as their centers,
(ii) for any two vertices $u,v\in N(R,S)$, $deg(u,v,S)\ge 3|R|$,  and
(iii) for any three vertices $u,v, w\in N(N(R,S),S)$, $deg(u,v,w, S)\ge 4|R|$.
Then  there is a ladder spanning on $R$ and some other $7|R|-2$ vertices from $S$
with $3|R|-1$ of them taken from $A$.
Particularly, the ladder has its first and last rungs in $S$.
\end{LEM}

The following simple observation is heavily used in the proofs
explicitly or implicitly.

\begin{LEM}\label{common-vertex}
Let $U=\{u_1,u_2\cdots, u_k\},S\subseteq V(G)$
be subsets. Then
$deg(u_1,u_2,\cdots, u_k, S)\ge |S|-(deg_{\overline{G}}(u_1,S)+ \cdots +deg_{\overline{G}}(u_k,S))
\ge |S|-k(|S|-\delta(U,S))=k\delta(U,S)-(k-1)|S|$.
\end{LEM}

Extremal Case 1 is easier than the other cases, so we start with it.

\subsection{Proof of Theorem~\ref{extremal1}}

We assume that $G$ has an approximate  vertex-cut $W$ with parameter $\beta$
such that $|W|\le 5\beta n$. Let $V_1$ and $V_2$ be the partition of $V-W$
such that $\delta[G[V_i]]\ge (1/2-6\beta)n$.
As $\delta(G)\ge (n+1)/2$, $(1/2-6\beta )n \le |V_i|\le (1/2+6\beta)n$.   We partition $W$ into two subsets as follows:
$$
W_1=\{w\in W\,|\, deg(w,V_1)\ge (n+1)/4-2.5\beta n\}\quad \mbox{and}\quad W_2=W-W_1.
$$
As $\delta(G)\ge (n+1)/2$, we have $deg(w,V_2)\ge (n+1)/4-2.5\beta n$ for any $w\in W_2$.
Since $G$ is 3-connected and $(1/2-6\beta)n>3$, there are three independent edges $p_1p_2$, $q_1q_2$,
and $r_1r_2$ between $G[V_1\cup W_1]$ and $G[V_2\cup W_2]$ with $p_1, q_1, r_1\in V_1\cup W_1$ and $p_2,q_2,r_2\in V_2\cup W_2$.

For $i=1,2$, by the partition of $W_i$, we see that $\delta(W_i, V_i)\ge 3|W_i|+3$.
Thus, $\delta(W_i, V_i-\{p_i, q_i\})\ge 3|W_i|$. There are
$|W_i-\{p_i,q_i\}|$ vertex-disjoint 3-stars with their centers in
$W_i-\{p_i,q_i\}$.
By Lemma~\ref{common-vertex},
we have
\begin{eqnarray*}
  deg(u,v, V_i-\{p_i,q_i\})&\ge &2\delta(G[V_i])-|V_i|\ge  (1/2-18\beta)n\ge 6|W_i|,\quad \mbox{for any $u,v\in V_i$;} \\
  deg(u,v, w, V_i-\{p_i,q_i\})&\ge & 3\delta(G[V_i])-2|V_i| \ge (1/2-30\beta)n\ge 7|W_i|, \quad \mbox{for any $u,v,w\in V_i$.}
\end{eqnarray*}
By Lemma~\ref{absorbing}, we can find a ladder $L_i$ which spans $W_i-\{p_i,q_i\}$
and another $7|W_i-\{p_i,q_i\}|-2$ vertices from $V_i-\{p_i, q_i\}$,  if $W_i-\{p_i, q_i\}\ne \emptyset$.
Denote $a_ib_i$ and $c_id_i$ the first and last rung of $L_i$\,(if $L_i$ exists),
respectively. Let
$$
G_i=G[V_i-V(L_i)] \quad \mbox{and}\quad n_i=|V(G_i)|.
$$
Then for $i=1,2$, if $x\in V(G_i)$ and $x\not\in \{p_i, q_i\}\cap W$,
\begin{equation*}\label{Gisize}
    n_i\ge (n+1)/2-6\beta n-7|W_i|\ge (n+1)/2-41\beta n,  \quad  deg_{G_i}(x)\ge \delta(G[V_i])-7|W_i|\ge (n+1)/2-41\beta n.
\end{equation*}
If $p_i\in W$, then $deg_{G_i}(p_i)\ge (n+1)/4-2.5\beta n-7|W_i|\ge (1/4-41\beta)n$. Similarly, if
 $q_i\in W$, then $deg_{G_i}(q_i)\ge (1/4-41\beta)n$.

Let $i=1,2$.
We now show that $G_i$ contains a spanning subgraph isomorphic to either $H_1$ or $H_2$ as defined in the previous section.
Since $n_i\le (1/2+6\beta)n$ and $deg_{G_i}(x)\ge (n+1)/2-41\beta n$ for any $x\in V(G_i)- W$,
any subgraph of $G_i$ induced by at least $(1/4-41\beta)n$
vertices  not in $W$ has minimum degree at least $(n+1)/2-41\beta n-(n_i-(1/4-41\beta)n)\ge(1/4-88\beta)n$, and thus has a matching of size at least 2.
Hence,  when $n_i$ is even, we can choose independent edges $e_i=x_iy_i$ and $f_i=z_iw_i$ with
$$
x_i, y_i \in \Gamma_{G_i}(p_i)-\{q_i\}\quad \mbox{and}\quad
z_i, w_i \in \Gamma_{G_i}(q_i)-\{p_i\}.
$$
And  if $n_i$ is odd,
we can choose independent edges $g_iy_i$\,(we may assume $g_i\ne r_i$),  $f_i=z_iw_i$, and a vertex $x_i$ with
$$
g_i, x_i, y_i \in \Gamma_{G_i}(p_i)-\{q_i\}, x_i\in \Gamma_{G_i}(g_i, y_i)-\{p_i,q_i\} \quad \mbox{and}\quad
z_i, w_i \in \Gamma_{G_i}(q_i)-\{x_i, p_i\},
$$
where the existence of the vertex $x_i$ is possible since  the subgraph of $G_i$ induced by $\Gamma_{G_i}(p_i)$
has minimum degree at least
$(1/2-41\beta)n-((1/2+6\beta)n-|\Gamma_{G_i}(p_i)|)\ge |\Gamma_{G_i}(p_i)|-47\beta n$,  and hence contains
a triangle.
In this case, again, denote  $e_i=x_iy_i$.   Let
$$
\begin{cases}
   G_i'=G_i-\{p_i,q_i\},   & \text{if $n_i$ is even}; \\
  G_i'=G_i-\{p_i, q_i,g_i\},    & \text{if $n_i$ is odd}.
\end{cases}
$$
By the definition above,  $|V(G_i')|$ is even.

The following claim is a modification of (1) of Lemma 2.2 in~\cite{MR2646098}.
\begin{CLA}\label{H12}
For $i=1,2$, let $a_{i}'b'_{i}, c_i'd_i'\in E(G_i')$ be two independent edges.
Then $G_i'$ contains two vertex disjoint ladders $Q_{i1}$ and $Q_{i2}$ spanning on $V(G_i')$
such that  $Q_{i1}$
has  $e_i=x_iy_i$  as its first rung,  $a_i'b_i' $ as its last rung,  and $Q_{i2}$ has $c_i'd_i'$ as its first rung and
 $f_i=z_iw_i$  as its last rung, where $e_i$ and $f_i$ are defined prior to this claim.
 \end{CLA}

\pf  We only show the claim for $i=1$ as the case for $i=2$ is similar.
Notice that by the definition of $G_1'$, $|V(G_1')|$ is even.
Since $|V(G_1')|\le (1/2+6\beta)n$ and $\delta(G_1')\ge (n+1)/2-41\beta n-3\ge |V(G_1')|/2+4$,
$G_1'$ has a perfect matching  $M$ containing $e_1, f_1, a_1'b_1', c_1'd_1'$.
We identify $a_1'$ and $c_1'$ into a vertex called  $s'$, and identify $b_1'$ and $d_1'$ into
a vertex called $t'$.
Denote $G_1''$ as the resulting graph. Note that
 $s't'\in E(G_1'')$ by the way of identifications.
Partition $V(G_1'')$
arbitrarily into $U_1$ and $U_2$ with $|U_1|=|U_2|$ such that $x_1, z_1, s'\in U_1$,
$y_1, w_1, t'\in U_2$, and let $M':=(M-\{a_1'b_1', c_1'd_1'\})\cup \{s't'\}\subseteq E_{G_1''}(U_1,U_2)$.
 Define an auxiliary graph $H'$ with vertex set $M'$ and edge set defined as follows.
If $xy, uv\in M'-\{s't'\}$ with $x,u\in U_1$ then $xy\sim_{H'} uv$ if and only if $x\sim_{G_1'} v$
and $y\sim_{G_1'} u$\,(we do not include the case that $x\sim_{G_1'} u$
and $y\sim_{G_1'} v$ as we defined a bipartition here). Particularly, for any $pq\in M'-\{s't'\}$ with $p\in U_1$,
 $pq\sim_{H'} s't'$  if and only if $p\sim_{G_1'} b_1', d_1'$
 and $q\sim_{G_1'} a_1', c_1'$.
Notice that  a ladder with rungs in $M'$ is corresponding to a path in $H'$ and vice versa.
Since $(1/2-41\beta)n-3\le |V(G_1')|\le (1/2+6\beta)n-2$ and $\delta(G_1')\ge (n+1)/2-41\beta n-3$,
any two vertices in $G_1'$ has at least $(1/2-88\beta) n-4\ge (1/2-89\beta) n$ common neighbors by Lemma~\ref{common-vertex}.
This together with the fact that $|U_1|=|U_2|\le |V(G_1'')|/2\le (1/4+3\beta)n $
gives that $\delta(U_1,U_2), \delta(U_2,U_1)\ge (1/4-92\beta)n$.
For each edge $uv\in M'$ with $u\in U_1$, $u$ is adjacent to at least $(1/4-92\beta)n$
other vertices in $U_2$ saturated by $M'$. Thus there are at least $(1/4-92\beta)n$ edges $\{u_jv_j\,|\,  v_j\in \Gamma_G(u,U_2)\}\subseteq M'$.
Among these vertices $\{u_j\,|\, v_j\in \Gamma_G(u,U_2)\}$ in $U_1$, at least $(1/4-92\beta)n-\left((1/4+3\beta)n-(1/4-92\beta)n\right)=(1/4-187\beta)n$
of them are neighbors of $v$. Thus, in $H'$, $uv$ is adjacent to at least $(1/4-187\beta)n$ neighbors, and thus
$$
\delta(H')\ge  (1/4-187\beta)n\ge |V(H')|/2+1,
$$
since $\beta <1/2200$ and $n$ is very large.
Hence $H'$ has a hamiltonian path starting with $e_1$, ending with $f_1$,
and having $s't'$  as an  internal vertex.
The path with $s't'$ replaced by $a_1'b_1'$ and $c_1'd_1'$  is corresponding to the required ladders in $G_1'$.
\qed

\begin{figure}[!htb]
\psfrag{x_1}{$x_1$} \psfrag{x_2}{$x_2$}
\psfrag{y_1}{$y_1$} \psfrag{y_2}{$y_2$} \psfrag{a_2}{$a_2$}\psfrag{a_2'}{$a_2'$}
\psfrag{b_2}{$b_2$}
\psfrag{c_2}{$c_2$}
\psfrag{d_2}{$d_2$}
\psfrag{a_1}{$a_1$}
\psfrag{b_1}{$b_1$}
\psfrag{c_1}{$c_1$}
\psfrag{d_1}{$d_1$}
\psfrag{L_2}{$L_2$}
\psfrag{L_1}{$L_1$}
\psfrag{b_2'}{$b_2'$}
\psfrag{c_2'}{$c_2'$}
\psfrag{d_2'}{$d_2'$}
\psfrag{a_1'}{$a_1'$}
\psfrag{b_1'}{$b_1'$}
\psfrag{c_1'}{$c_1'$}
\psfrag{d_1'}{$d_1'$}\psfrag{d}{$d$} \psfrag{z_1}{$z_1$}
\psfrag{z_2}{$z_2$}
\psfrag{w_1}{$w_1$}
\psfrag{w_2}{$w_2$}
\psfrag{g_2}{$g_2$}
\psfrag{p_1}{$p_1$}
\psfrag{p_2}{$p_2$}
\psfrag{p_3}{$p_3$}
\psfrag{q_1}{$q_1$}
\psfrag{q_2}{$q_2$}
\psfrag{r_1}{$r_1$}
\psfrag{r_2}{$r_2$}
\psfrag{s_1}{$s_1$}
\psfrag{s_2}{$s_2$}
\psfrag{q_3}{$q_3$}
\psfrag{Q_11'}{$Q_{11}$}\psfrag{Q_12'}{$Q_{12}$}\psfrag{Q_21'}{$Q_{21}$}\psfrag{Q_22'}{$Q_{22}$}
\psfrag{L^T_1}{$L^T_1$}
\psfrag{F}{$F$}
\begin{center}
  \includegraphics[scale=0.4]{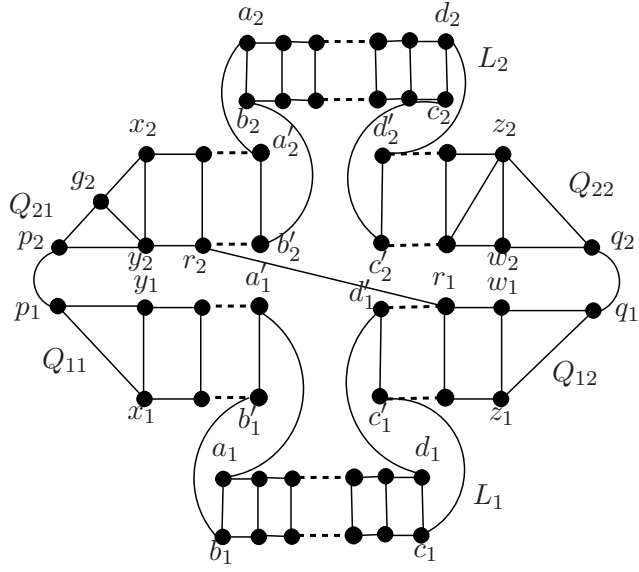}\\
\end{center}
\vspace{-3mm}
  \caption{A Halin graph $H$ }\label{halin2}
\end{figure}

We may assume $n_1$ is even and $n_2$ is odd and construct a spanning Halin subgraph of
$G$\,(the construction for the other three cases follow a similar argument).
Recall that $p_1p_2, q_1q_2, r_1r_2$ are the three prescribed independent edges
between $G[V_1\cup W_1]$ and $G[V_2\cup W_2]$, where
 $p_1,q_1, r_1\in V_1\cup W_1$ and $p_2, q_2, g_2, r_2\in V_2\cup W_2$.
For a uniform discussion, we may assume that both of the ladders $L_1$
and $L_2$ exist. Let $i=1,2$. Recall that $L_i$ has $a_ib_i$ as its first rung and $c_id_i$
as its last rung.
Choose $a_i'\in \Gamma_{G}(a_i, V(G_i'))$, $b_i'\in \Gamma_{G}(b_i, V(G_i'))$
such that $a_i'b_i'\in E(G)$ and
$c_i'\in \Gamma_{G}(c_i, V(G_i'))$, $d_i'\in \Gamma_{G}(d_i, V(G_i'))$
such that $c_i'd_i'\in E(G)$\,
($a_i', b_i', c_i', d_i'$ are chosen mutually distinct and distinct from
$x_i, y_i, z_i, w_i, g_i, r_i$). This is possible as $\delta(V_i, V(G_i'))\ge (n+1)/2-41\beta n-2$.
Let $Q_{i1}$ and $Q_{i2}$ be the ladders of $G_i'$ given by Claim~\ref{H12}.
Set $H_a=Q_{11}L_1 Q_{12}\cup \{p_1x_1, p_1y_1, q_1z_1, q_1w_1\}$.
Assume $Q_{21}L_2 Q_{22}$ is a ladder can be denoted as
$\overrightarrow{x_2y_2}-Q_{21}L_2 Q_{22}-\overrightarrow{z_2w_2}$.
To make $r_2$ a Halin constructible vertex, we let
$H_b=Q_{21}L_2 Q_{22}\cup \{g_2x_2,  g_2y_2, p_2g_2, p_2y_2, q_2z_2, q_2w_2\}$
if $r_2$ is on the shortest $(y_2,w_2)$-path in $Q_{21}L_2 Q_{22}$,
and let $H_b=Q_{21}L_2 Q_{22}\cup \{g_2x_2,  g_2y_2, p_2g_2, p_2x_2, q_2z_2, q_2w_2\}$
if $r_2$ is on the shortest $(x_2,z_2)$-path\,(recall that $g_2,x_2,y_2\in \Gamma_{G_2}(p_2)$).
Let $H=H_a\cup H_b\cup \{p_1p_2, r_1r_2, q_1q_2\}$.
Then $H$
is a spanning Halin subgraph of $G$ by Proposition~\ref{Prop:Halin_H1_H2} as
$H_a\cup p_1q_1\cong H_1$ and $H_b\cup p_2q_2\cong H_2$.
 Figure~\ref{halin2} gives a construction of $H$ for the above case when
 $r_2$ is on the shortest $(y_2,w_2)$-path in $Q_{21}L_2 Q_{22}$.

\subsection{Proof of Theorem~\ref{extremal2}}

Recall Extremal Case 2: There exists a partition $V_1\cup V_2$ of $V$
such that $|V_1|\ge (1/2-7\beta)n$ and $\Delta(G[V_1])\le \beta n$.
Since $\delta(G)\ge (n+1)/2$,  the assumptions imply that
$$
 (1/2-7\beta)n \le |V_1| \le (1/2+\beta)n  \quad \mbox{and} \quad
 (1/2-\beta)n \le |V_2| \le (1/2+7\beta)n.
$$
Let  $\beta$ and  $\alpha$ be real numbers satisfying  $\beta \le \alpha/20$ and $ \alpha\le (1/17)^3$.
Set $\alpha_1=\alpha^{1/3}$ and  $\alpha_2=\alpha^{2/3}$.
We first repartition
$V(G)$ as follows.
\begin{eqnarray*}
V_2'& = & \{v\in V_2\,|\, deg(v,V_1)\ge (1-\alpha_1)|V_1|\} ,  V_{01}=\{v\in V_2-V_2'\,|\, deg(v,V_2')\ge (1-\alpha_1)|V_2'|\},  \\
V_1' & = & V_1\cup V_{01}, \quad  \mbox{and}\quad V_0=V_2-V_2'-V_{01}.
\end{eqnarray*}
\begin{CLA}\label{V2-V2'_size}
$|V_{01}\cup V_0|=|V_2-V_2'|\le \alpha_2|V_2|$.
\end{CLA}

\pf
Notice that $e(V_1, V_2)\ge (1/2-\beta)n |V_1| \ge \frac{1/2-\beta}{1/2+7\beta}|V_1||V_2|\ge (1-\alpha)|V_1||V_2|$ as $\beta \le \alpha/20$. Hence,
\begin{eqnarray*}
(1-\alpha)|V_1||V_2|\le & e(V_1, V_2) & \le e(V_1, V_2')+e(V_1, V_2-V_2')\le |V_1||V_2'|+(1-\alpha_1)|V_1||V_2-V_2'|.
\end{eqnarray*}
This gives that $|V_{01}\cup V_0|=|V_2-V_2'|\le \alpha_2|V_2|$.
\qed

As a result of moving vertices from $V_2$ to $V_1$ and
by Claim~\ref{V2-V2'_size}, we have the following.
\begin{eqnarray}
(1/2-7\beta)n\le |V_1'|&\le & (1/2+\beta)n+|V_{01}|\le (1/2+\beta)n+\alpha_2(1/2+7\beta)n\le (1/2+\alpha_2)n , \nonumber\\
(1/2-\alpha_2)n\le |V_2'|&\le & (1/2+7\beta)n, \nonumber\\
\delta(V_1', V_2') & \ge  & \min\{(1/2-\beta)n-|V_2-V_2'|, (1-\alpha_1)|V_2'|\}\ge (1/2-2\alpha_1/3)n , \nonumber\\
\delta(V_2', V_1') &\ge  &(1-\alpha_1)|V_1|\ge (1-\alpha_1)(1/2-7\beta)n\ge (1/2-2\alpha_1/3)n, \label{degrees} \\
\delta(V_{0}, V_1')&\ge & (n+1)/2-(1-\alpha_1)|V_2'|-|V_0| \ge  \alpha_1 n/3\ge 6|V_0|+20, \nonumber \\
\delta(V_{0}, V_2')&\ge & (n+1)/2-(1-\alpha_1)|V_1|-|V_0\cup V_{01}| \ge   \alpha_1 n/3\ge 6|V_0|+20.  \nonumber
\end{eqnarray}

\begin{CLA}\label{nowheel}
We may assume that $\Delta(G)<n-1$.
\end{CLA}

\pf Suppose on the contrary and let $w\in V(G)$ such that $deg(w)=n-1$. Then by $\delta(G)\ge (n+1)/2$
we have $\delta(G-w)\ge (n-1)/2$, and thus $G-w$ has a hamiltonian cycle. This implies that $G$ has
a spanning wheel subgraph,  in particular, a spanning Halin subgraph of $G$.
\qed

\begin{CLA}\label{subgraph_K}
There exists a subgraph $T\subseteq G$ with $|V(T)|\equiv n\,(\mbox{mod}\,\,2)$
such that $T$ and $G-V(T)$ satisfy the following conditions.
\begin{itemize}
 \item[$($i$)$] $T$ is isomorphic to some graph in $\{T_1,T_2,\cdots, T_5\}$;
  \item [$($ii$)$] Let $2m=n-|V(T)|$. Then $G-V(T)$ contains a balanced spanning bipartite graph $G'$
with partite sets $U_1$ and $U_2$ such that $|U_1|=|U_2|=m$.
\item [$($iii$)$]
There exists a subset $W$ of $U_1\cup U_2$ with at most $\alpha_2n$ vertices
such that $deg_{G'}(x, V(G')-W)\ge (1-\alpha_1-2\alpha_2)m$ for all $x\not\in W$.
  \item[$($iv$)$] Assume that $T$ has head link $x_1x_2$ and tail link $y_1y_2$.
  There exist $x_1'x_2', y_1'y_2'\in E(G')$ such that
  $x_i', y_i'\in U_i-W$,   $x_{3-i}'\sim x_i$,  and $y_{3-i}'\sim y_i$, for $i=1,2$;
  and if $T$ has a pendent vertex, then the pendent vertex is contained in $V'_1\cup V'_2-W$.
   \item[$($v$)$]There are $|W|$ vertex-disjoint 3-stars  in $G'-\{x_1',x_2', y_1', y_2'\}$
  with the vertices in $W$ as their centers.
  \end{itemize}
\end{CLA}

\pf  By~\eqref{degrees} and Lemma~\ref{common-vertex},  for $i=1,2$, we notice that for any $u,v,w\in V_i'$,
\begin{eqnarray}\label{common}
deg(u,v,w, V_{3-i}') & \ge & 3\delta(V_i', V_{3-i}')-2|V_{3-i}'|\ge (1/2-3\alpha_1)n>n/4.
\end{eqnarray}
We now separate the proof into two cases according to the parity of  $n$.

{\noindent \bf Case 1. $n$ is even. }

Suppose first that $\max\{|V_1'|, |V_2'|\}\le n/2$.
We arbitrarily partition $V_0$
into $V_{10}$ and $V_{20}$ such that $|V'_1\cup V_{10}|=|V'_2\cup V_{20}|=n/2$.
Suppose $G[V_1']$ contains an edge $x_1u_1$
and there is a  vertex $u_2\in \Gamma(u_1, V_2')$ such that $u_2$ is adjacent to a vertex $y_2\in V_2'$.
By~\eqref{common},  there exist distinct vertices  $x_2\in \Gamma(x_1, u_1, V_2')-\{y_2,u_1\}, y_1\in \Gamma(y_2, u_2, V_1')-\{x_1,u_1\}$.
Then $G[\{x_1, u_1, x_2, y_1, u_1, y_2\}]$ contains a subgraph $T$ isomorphic to $T_1$.
So we assume $G[V_1']$ contains an edge $x_1u_1$ and
no vertex in  $\Gamma(u_1, V_2')$ is adjacent to any vertex in $V_2'$.
As $\delta(G)\ge (n+1)/2$,  $\delta(G[V'_2\cup V_{20}])\ge 1$.
Let $u_2\in \Gamma(u_1, V_2') $ and $u_2y_2\in E(G[V'_2\cup V_{20}])$.
Since $deg(u_2, V_1')\ge (n+1)/2-|V_0|>|V_1'\cup V_{10}|-|V_0|$ and $deg(y_2, V_1')\ge 3|V_0|+10$,
$\deg(u_2, y_2, V_1'\cup V_{10})\ge 2|V_0|+10$.
Let $x_2\in \Gamma(x_1, u_1, V_2')-\{y_2,u_2\}, y_1\in \Gamma(y_2, u_2, V_1')-\{x_1,u_1\}$.
Then $G[\{x_1, u_1, x_2, y_1, u_2, y_2\}]$ contains a subgraph $T$ isomorphic to $T_1$.
By symmetry, we can find $T\cong T_1$ if
$G[V_2']$ contains an edge.  Hence we assume that both
$V_1'$ and $V_2'$ are independent sets.
Again, as $\delta(G)\ge (n+1)/2$,  $\delta(G[V'_1\cup V_{10}]), \delta(G[V'_2\cup V_{20}])\ge 1$.
Let $x_1u_1\in E(G[V'_1\cup V_{10}])$ and $y_2u_2\in E(G[V'_2\cup V_{20}])$
such that $x_1\in V_1'$ and $u_2\in \Gamma(u_1, V_2')$.
Since $deg(x_1, V_2')\ge (n+1)/2-|V_0|>|V_2'\cup V_{20}|-|V_0|$ and $deg(u_1,V_2')\ge 3|V_0|+10$,
we have $deg(x_1,u_1, V_2')\ge 2|V_0|+10$.
Hence, there exists  $x_2\in \Gamma(x_1, u_1, V_2')-\{y_2,u_2\}$.
Similarly,
there exists  $ y_1\in \Gamma(y_2, u_2, V_1')-\{x_1,u_1\}$.
Then $G[\{x_1, u_1, x_2, y_1, u_2, y_2\}]$ contains a subgraph $T$ isomorphic to $T_1$.
Let $m=(n-6)/2, U_1=(V_1'\cup V_{10})-V(T)$ and $U_2=(V_2'\cup V_{20})-V(T)$,
and $W=V_0-V(T)$.  We then have $|U_1|=|U_2|=m$.

Let $G'=(V(G)-V(T), E_G(U_1,  U_2))$ be the bipartite graph with partite sets $U_1$ and $U_2$.
Notice that  $|W|\le |V_0|\le \alpha_2|V_2|<\alpha_2n$.
By~(\ref{degrees}),  we have $deg_{G'}(x, V(G')-W)\ge (1-\alpha_1-2\alpha_2)m$ for all $x\notin W$.
This shows (iii). By the construction of $T$ above, we have  $x_1, y_1\in V_1'-W$. Let $i=1, 2$.
By~(\ref{degrees}), we have $\delta(V_0, U_i-W)=\delta(V_0, V_i'-V(T))\ge 3|V_0|+6$.
So $|\Gamma_{G'}(y_2, U_1-W)|, |\Gamma_{G'}(x_2, U_1-W)|\ge 3|V_0|+6$.
 Applying statement (iii) and Lemma~\ref{common-vertex}, we have $e_{G'}(\Gamma_{G'}(x_2, U_1-W),\Gamma_{G'}(x_1, U_2-W)),
 e_{G'}(\Gamma_{G'}(y_2, U_1-W), \Gamma_{G'}(y_1, U_2-W))\ge (3|V_0|+6)(1-2\alpha_1-4\alpha_2)m>2m$.
 Hence, we can find  independent edges $x_1'x_2'$ and $y_1'y_2'$ such that
  $x_i', y_i'\in U_i-W$,   $x_{3-i}'\sim x_i$,  and $y_{3-i}'\sim y_i$. This gives statement (iv).
  Finally, as  $\delta(V_0, U_i-W)\ge 3|V_0|+6$, we have
  $\delta(V_0, U_i-W-\{x_1',x_2', y_1', y_2'\})\ge 3|V_0|+2$.
 Hence,  there are  vertex-disjoint 3-stars with their centers in $W$.

Otherwise we have $\max\{|V_1'|, |V_2'|\}> n/2$.
By~(\ref{degrees}), we have the same lower bound for
$\delta(V_1',V_2')$, $\delta(V_2',V_1')$,
and $\delta(V_0,V_1')$, $\delta(V_0,V_2')$.
Furthermore, all the argument in the following will depend
only on the degree conditions, so we assume,
 w.l.o.g., that $|V_1'|\ge n/2+1$.
Then $\delta(G[V_1'])\ge 2$ and thus $G[V_1']$
 contains two vertex-disjoint paths isomorphic to $P_3$ and $P_2$, respectively.
Let $m=(n-8)/2$.  We consider three cases here.
Case (a): $|V_1'|-5\le m$.  Then
 let $x_1u_1w_1, y_1v_1\subseteq G[V_1']$
be  two vertex-disjoint paths, and let
$x_2\in \Gamma(x_1, u_1, w_1, V_2'), y_2\in \Gamma(y_1, v_1, V_2')$
and $z\in \Gamma(w_1, v_1, V_2')$ be three distinct vertices.
Then $G[\{x_1, u_1, w_1, x_2, z, y_1, v_1, y_2\}]$
contains a subgraph $T$ isomorphic to $T_4$.  Notice that $|V_2'-V(T)|\le m$.
We
arbitrarily partition $V_0$
into $V_{10}$ and $V_{20}$ such that $|V'_1\cup V_{10}|=|V'_2\cup V_{20}|=m$.
Let $U_1=(V_1' \cup V_{10})-V(T)$, $U_2=(V_2'\cup V_{20})-V(T)$, and $W=V_0$.
Hence we assume $|V_1'|-5 = m+t_1$
for some $t_1\ge 1$. This implies that $|V_1'|= n/2+t_1+1$ and thus $\delta(G[V_1'])\ge t_1+2$.
Let $V_1^0$ be the set of vertices $u \in V_1'$
such that $deg(u, V_1')\ge \alpha_1 m$.
Case (b): $|V_1^0|\ge |V_1'|-5-m$.  Then we form a set $W$ with $|V'_1|-5-m$
vertices from $V_1^0$ and all the vertices of $V_0$.
Then $|V_1'-W|=m+5+t_1-(|V_1'|-5-m)=m+5=n/2+1$, and hence
$\delta(G[V_1'-W])\ge 2$.  Similarly as in Case (a), we can find
a subgraph $T$ of $G$ contained in  $G[(V_1'\cup V_2')-W]$ isomorphic to $T_4$.
Let $U_1=V_1'-V(T)-W$, $U_2=(V_2'\cup W)-V(T)$. Then $|U_1|=|U_2|=m$.
Thus we have Case (c): $|V_1^0|< |V_1'|-5-m$.
Suppose that $|V_1'-V_1^0|=m+5+t_1'=n/2+t_1'+1$ for some $t_1'\ge 1$.
This implies that $\delta(G[V_1'-V_1^0])\ge t_1'+2$.

We show that
$G[V_1'-V_1^0]$ contains $t_1'+2$ vertex-disjoint 3-stars. To see this,
suppose $G[V_1'-V_1^0]$ contains a subgraph $M$ of at most $s<t_1'+2$
3-stars. By counting the number of edges between $V(M)$ and $V_1'-V_1^0-V(M)$
in two ways, we get that $t_1'|V_1'-V_1^0-V(M)|\le e_{G-V_1^0}(V(M), V_1'-V_1^0-V(M))\le 4s \Delta(G[V_1'-V_1^0])\le 4s \alpha_1m$.
Since $|V_1'-V_1^0|=m+5+t_1'=n/2+t_1'+1$, $|V_1'-V_1^0-V(M)|\ge m -3t_1'\ge m-6\alpha_2m$, where the last inequality holds
as $|V_1'|\le (1/2+\beta)n +\alpha_2|V_2'|$ implying that
$t_1'\le |V_1'|-m -5\le 2\alpha_2m$.  This, together with the assumption that $\alpha\le (1/8)^3$ gives that
$s\ge t_1'+2$, showing a contradiction. Hence we have $s\ge t_1'+2$. Let
$x_1u_1w_1$ and $y_1v_1$ be two paths taken from two 3-stars in $M$. Then
we can find  a subgraph $T$ of $G$  isomorphic to $T_4$ in the same way as in Case (a).
We take exactly $t_1'$ 3-stars from the remaining ones in $M$ and denote the centers
of these stars by $W'$.
Let $U_1=V_1'-V_1^0-V(T)-W'$,  $W=W'\cup V_1^0\cup V_0$, and $U_2=(V_2'\cup W)-V(T)$. Then $|U_1|=|U_2|=m$.

For the partition of $U_1$ and $U_2$ in all the cases discussed in the paragraph above, we let
$G'=(V(G)-V(T), E_G(U_1,  U_2))$ be the bipartite graph with partite sets $U_1$ and $U_2$.
Notice that $|W|\le |V_0|\le \alpha_2 n$ if Case (a) occurs,
$|W|\le |V_0|+|V_1'|-m-5\le (1/2+\beta)n+|V_0\cup V_{01}|-n/2-1\le \alpha_2n$ if Case (b) occurs, and
$|W|=|W'\cup V_1^0\cup V_0|=|V_1'-U_1-V(T)|+|V_0\cup V_{01}|\le (1/2+\beta)n-(1/2-4)n+|V_0\cup V_{01}|\le \alpha_2n$
if Case (c) occurs. (Recall that $|V_1'|\le (1/2+\beta)n+|V_{01}|$ and $|V_0\cup V_{01}|\le \alpha_2|V_2|$ from~\eqref{degrees}.)
Since $\delta(V_2', V_1') \ge (1-2\alpha_1/3)n$ from~\eqref{degrees} and $|V_1'-U_1|\le 2\alpha_2m$,
we have $\delta(U_2-W, U_1-W)\ge (1-\alpha_1-2\alpha_2)m$. On the other hand,
from~\eqref{degrees}, $\delta(V_1', V_2') \ge (1/2-2\alpha_1/3)n$. This gives that
$\delta(U_1-W, U_2-W)\ge (1-\alpha_1-2\alpha_2)m$. Hence, we have
$deg_{G'}(x, V(G')-W)\ge (1-\alpha_1-2\alpha_2)m$ for all $x\notin W$.
According to the construction of $T$, we have $x_1,y_1\in V_1'-W$.
Applying statement (iii)  and Lemma~\ref{common-vertex}, we have $e_{G'}(\Gamma_{G'}(x_1, U_2-W), \Gamma_{G'}(x_2, U_1-W)),
 e_{G'}(\Gamma_{G'}(y_1, U_2-W), \Gamma_{G'}(y_2, U_1-W))\ge (3|V_0|+6)(1-2\alpha_1-4\alpha_2)m>2m$.
 Hence, we can find  independent edges $x_1'x_2'$ and $y_1'y_2'$ such that
 $x_i', y_i'\in U_i-W$,   $x_{3-i}'\sim x_i$,  and $y_{3-i}'\sim y_i$. By the construction of
 $T$, $T$ is isomorphic to $T_4$, and the pendent vertex $z\in V_2'\subseteq V'_1\cup V_2'-W$.  This gives statement (iv).
 Finally, as
 \begin{eqnarray*}
   \delta(V_0, U_1-W) &\ge &  \delta(V_0, V_1')-|V_1'-(U_1-W)| \ge \alpha_1n/3-(1/2+\alpha_2)n+n/2-4-\alpha_2n \\
    &\ge & (1/3\alpha_1-2\alpha_2)n-4\ge 3|W|+5,
 \end{eqnarray*}
 we have
 $\delta(V_0, U_1-W-\{x_1', x_2', y_1', y_2'\})\ge 3|W|+1$.
By the definition of  $V_1^0$, we have $\delta(V_1^0, V_1'-W-\{x_1', x_2', y_1', y_2'\})\ge \alpha_1m-\alpha_2 n-4\ge 3|W|$.
For the vertices in $W'$ in Case (c), we already know that there are vertex-disjoint
3-stars in $G'$ with centers in $W'$.
 Hence,  regardless of the construction of $W$, we can
 always find  vertex-disjoint 3-stars with their centers in $W$.


{\noindent \bf Case 2. $n$ is odd. }

Suppose first that $\max\{|V_1'|, |V_2'|\}\le (n+1)/2$
and let $m=(n-7)/2$.
We arbitrarily partition $V_0$
into $V_{10}$ and $V_{20}$ such that, w.l.o.g.,  say  $|V'_1\cup V_{10}|=(n+1)/2$ and $|V'_2\cup V_{20}|=(n-1)/2$.
(Again, here we use the symmetry of the lower bounds on
$\delta(V_1',V_2')$, $\delta(V_2',V_1')$,
and $\delta(V_0,V_1')$, $\delta(V_0,V_2')$ from~(\ref{degrees}).)
We show that $G[V_1'\cup V_{10}]$ either contains two independent edges or is isomorphic to $K_{1,(n-1)/2}$.
As $\delta(G)\ge (n+1)/2$, we have $\delta(G[V_1'\cup V_{10}])\ge 1$.
Since $n$ is sufficiently large, $(n+1)/2>3$.  Then it is easy to see that
if $G[V_1'\cup V_{10}]\not\cong K_{1,(n-1)/2}$, then $G[V_1'\cup V_{10}]$ contains
two independent edges. Furthermore, we can choose two independent edges
$x_1u_1$  and $y_1v_1$ such that $u_1, v_1\in V_1'$. This is obvious if
$|V_{10}|\le 1$. So we assume $|V_{10}|\ge 2$. As $\delta(V_0, V_1')\ge 3|V_0|+10$,
by choosing $x_1,y_1\in V_{10}$, we can choose distinct vertices $u_1\in \Gamma(x_1, V_1')$ and $v_1\in \Gamma(y_1, V_1')$.
Let $x_2\in \Gamma(x_1, u_1, V_2'),  y_2\in \Gamma(y_1, v_1, V_2')$
and $z\in \Gamma(u_1, v_1, V_2')$.
Then $G[\{x_1, u_1, x_2, y_1, v_1, y_2,z \}]$ contains a subgraph $T$ isomorphic to $T_3$.
We assume now that  $G[V_1'\cup V_{10}]$
is isomorphic to $K_{1, (n-1)/2}$.
Let $u_1$ be the center of the star $K_{1, (n-1)/2}$.
Then each leaf of the star has at least $(n-1)/2$
neighbors in $V_2'\cup V_{20}$. Since $|V_2'\cup V_{20}|=(n-1)/2$,
we have $\Gamma(v, V_2'\cup V_{20})=V_2'\cup V_{20}$
if $v\in V_1'\cup V_{10}-\{ u_1\}$.
By the definition of
$V_0$,  $\Delta(V_0, V_{1}')< (1-\alpha_1)|V_1|+|V_{01}|$ and $\Delta(V_0, V_2')<(1-\alpha_1)|V_2'|$, and so  $u_1\in V_1'$, $V_{10}=\emptyset$ and $V_{20}=\emptyset$.  We claim that $V_2'$ is not an independent
set. Otherwise, by $\delta(G)\ge (n+1)/2$, for each $v\in V_2'$, $\Gamma(v, V_1')=V_1'$.
This in turn shows that $u_1$ has degree $n-1$, showing a contradiction to
Claim~\ref{nowheel}.
So let  $y_2v_2\in E(G[V_2'])$ be an edge.
Let  $w_1\in \Gamma(v_2, V_1')-\{u_1\}$
 and  $w_1u_1x_1$ be the path containing
 $w_1$.  Choose $y_1\in \Gamma(y_2, v_2, V_1')-\{w_1,u_1,x_1\}$
 and $x_2\in \Gamma(x_1, u_1, w_1, V_2')-\{y_2,v_2\} $.
 Then $G[\{x_1, u_1, x_2, w_1, v_2, y_2, y_1\}]$ contains a subgraph $T$ isomorphic
 to  $T_2$.
 Let $U_1=(V_1'\cup V_{10})-V(T)$ and $U_2=(V_2'\cup V_{20})-V(T)$ and $W=V_0-V(T)$.
 We have $|U_1|=|U_2|=m$ and $|W|\le |V_0|\le \alpha_2n$.


Otherwise we have $\max\{|V_1'|, |V_2'|\}\ge (n+1)/2+1$.
By the symmetry of lower bounds on
degrees  related to $V_1'$ and $V_2'$ from~(\ref{degrees}),
 we assume, w.l.o.g., that $|V_1'|\ge (n+1)/2+1$.
Then $\delta(G[V_1'])\ge 2$ and thus $G[V_1']$  contains two independent edges.
Let $m=(n-7)/2$
and $V_1^0$ be the set of vertices $u \in V_1'$
such that $deg(u, V_1')\ge \alpha_1 m$.
 Since $|V_1'|\ge (n+1)/2+1>m+4$,
we assume $|V_1'|=m+4+t_1$ for some $t_1\ge 1$.
We consider three  cases here.

Case (a): $|V_1^0|\ge |V_1'|-m-4$. \quad
We form a set $W$ with $|V'_1|-4-m$
vertices from $V_1^0$ and all the vertices of $V_0$.
Then $|V_1'-W|=m+4+t_1-(|V_1'|-4-m)=m+4=(n+1)/2$,
and we have $\delta(G[V_1'-W])\ge 1$.  As any vertex
 $u\in V_1'-W$ is a vertex such that
 $deg(u, V_1')< \alpha_1 m$,
we know $G[V_1'-W]$
contains two independent edges.
Let $x_1u_1, y_1v_1\subseteq E(G[V_1'-W])$
be two independent edges, and
let  $x_2\in \Gamma(x_1, u_1,  V_2'), y_2\in \Gamma(y_1, v_1, V_2')$
and $z\in \Gamma(w_1, v_1, V_2')$ be three distinct vertices.
Then $G[\{x_1, u_1,  x_2, z, y_1, v_1, y_2\}]$ contains a subgraph $T$ isomorphic to $T_3$.
Let $U_1=V_1'-V(T)-W$, $U_2=(V_2'\cup W)-V(T)$. Then $|U_1|=|U_2|=m$
and $|W|\le |V_0|+|V_1'-U_1|\le |V_2-V_2'|+\beta n +4\le \alpha_2n$.

Thus we have $|V_1^0|< |V_1'|-4-m$.
Suppose that $|V_1'-V_1^0|=m+4+t_1'=(n+1)/2+t_1'$ for some $t_1'\ge 1$.
This implies that $\delta(G[V_1'-V_1^0])\ge t_1'+1$.

Case (b): $t_1'\ge 2$. \quad
We show that
$G[V_1'-V_1^0]$ contains $t_1'+2$ vertex-disjoint 3-stars.
To see this, suppose $G[V_1'-V_1^0]$ contains a subgraph  $M$ of
at most $s$ vertex disjoint 3-stars. We may assume that $s< t_1'+2$. Then we have
$(t_1-1) |V_1'-V_1^0-V(M)|\le e_{G-V_1^0}(V(M), V_1'-V_1^0-V(M))\le 4s\Delta(G[V_1'-V_1^0])$.
Since $|V_1'-V_1^0|=m+4+t_1'=(n+1)/2+t_1'$, $|V_1'-V_1^0-V(M)|\ge m -3t_1'\ge m-6\alpha_2m$, where the last inequality holds
as $|V_1'|\le (1/2+\beta)n +\alpha_2 |V_2'|$ implying that
$t_1'\le |V_1'|-m -4\le 2\alpha_2m$.  This, together with the assumption that $\alpha\le (1/8)^3$ gives that
$s\ge t_1'+2$, showing a contradiction. Hence we have $s\ge t_1'+2$.
Let
$x_1u_1$ and $y_1v_1$ be two paths taken from two 3-stars in $M$. Then
we can find  a subgraph $T$ of $G$  isomorphic to $T_3$ the same way
as in Case (a).
We take exactly $t_1'$ 3-stars from the remaining ones in $M$ and denote the centers
of these stars by $W'$.
Let $U_1=V_1'-V_1^0-V(T)-W'$,  $W=W'\cup V_1^0\cup V_0$, and $U_2=(V_2'\cup W)-V(T)$. Then $|U_1|=|U_2|=m$.

Case (c): $t_1'=1$.\quad   In this case, we let $m=(n-9)/2$.
If  $G[V_1'-V_1^0]$ contains a vertex adjacent to all other vertices in $V_1'-V_1^0$,
then the vertex would be contained in $V_1^0$ by the definition of $V_1^0$.
Hence,
we assume that $G[V_1'-V_1^0]$ has no vertex adjacent to all other vertices in $V_1'-V_1^0$.
Then by the assumptions that $\delta(G)\ge (n+1)/2$ and $|V_1'-V_1^0|=(n+1)/2+1$,
we can find two copies of vertex disjoint $P_3$s in $G[V_1'-V_1^0]$. Let
$x_1u_1w_1$ and $y_1v_1z_1$
be two $P_3$s in $G[V_1'-V_1^0]$.
There exist  distinct vertices $x_2\in \Gamma(x_1, u_1, w_1, V_2'), y_2\in \Gamma(y_1, v_1, z_1, V_2')$
and $z\in \Gamma(w_1,z_1, V_2')$.
Then $G[\{x_1, u_1, w_1, x_2,  y_1, v_1,z_1, y_2,z\}]$
contains a subgraph $T$ isomorphic to $T_5$.
Let $U_1=V_1'-V_1^0-V(T)$,  $W=V_1^0\cup V_0$, and $U_2=(V_2'\cup W)-V(T)$. Then $|U_1|=|U_2|=m$.

For the partition of $U_1$ and $U_2$ in all the cases discussed in Case 2, we let
$G'=(V(G)-V(T), E_G(U_1,  U_2))$ be the bipartite graph with partite sets $U_1$ and $U_2$.
Similarly as in Case 1, we can show that all the statements (i)-(v) hold.
  \qed

Let $W_1=U_1\cap W$ and $W_2=U_2\cap W$.
By (v) of Claim~\ref{subgraph_K}, we know
that
 there are $|W_1|$ vertex-disjoint 3-stars with centers in $W_1$
and all other vertices in  $U_2-W_2-\{x_1', y_1', x_2', y_2'\}$,
and $|W_2|$ vertex-disjoint 3-stars with centers in $W_2$
and all other vertices in  $U_1-W_1-\{x_1', y_1', x_2', y_2'\}$,
and all these $|W_1|+|W_2|$ stars are vertex-disjoint.
Let $S$ be the union of the 3-stars with centers in $W_2$.

For any $u,v\in U_2-W_2$, $\Gamma(u,v, U_1-W_1-V(S)-\{x'_1, x_2', y_1', y_2'\})\ge 3|W_1|$,
and for any $u,v, w\in U_1-W_1-V(S)$, $\Gamma(u,v, w, U_2-V(S)-\{x'_1, x_2', y_1', y_2'\})\ge 4|W_1|$.
By Lemma~\ref{absorbing2}, we can
find a ladder $L_1$ disjoint from the 3-stars in $S$ with centers in $W_2$
such that $L_1$ is
spanning on $W_1$,
$4|W_1|-1$ vertices from  $U_2-W_2-\{x'_1, x_2', y_1', y_2'\}$,
and another $3|W_1|-1$ vertices from  $U_1-W_1-\{x'_1, x_2', y_1', y_2'\}$,
 if $W_1\ne \emptyset$.

For any $u,v\in U_1-W_1$, $\Gamma(u,v, U_2-W_2-V(L_1)-\{x'_1, x_2', y_1', y_2'\})\ge 3|W_2|$,
and for any $u,v, w\in U_2-W_2$, $\Gamma(u,v, w, U_1-W_1-V(L_1)-\{x'_1, x_2', y_1', y_2'\})\ge 4|W_2|$.
By Lemma~\ref{absorbing2}, we can
find a ladder $L_2$ disjoint from $L_1$
such that $L_2$ is
spanning on $W_2$,
$4|W_2|-1$ vertices from  $U_1-V(L_1)-\{x'_1, x_2', y_1', y_2'\}$,
and another $3|W_2|-1$ vertices from  $U_2-V(L_1)-\{x'_1, x_2', y_1', y_2'\}$,
 if $W_2\ne \emptyset$.

%
%
Denote $a_{1i}a_{2i}$ and $b_{1i}b_{2i}$ the first and last rungs of $L_i$\,(if $L_i$ exists),
respectively, where $a_{1i}, b_{1i}\in U_1$.  As $|U_1|=|U_2|$, and we
took $4|W_1|+4|W_2|-2$ vertices respectively from $U_1$ and
$U_2$ when constructing $L_1$ and $L_2$, we have $|U_1-V(L_1\cup L_2)|=|U_2-V(L_1\cup L_2)|$.
Let
$$
U'_i=U_i- V(L_1\cup L_2), \quad \mbox{}\quad m'=|U_1'|=|U_2'|, \quad \mbox{and} \quad G''=G''(U_1'\cup U_2', E_G(U_1', U_2')).
$$
Since $|W|\le \alpha_2n$,  $m\ge (n-9)/2$,
and  $n$ is sufficiently large, we have $1/n+7|W|\le 15\alpha_2m$.
As $\delta(G'-W)\ge (1-\alpha_1-2\alpha_2)m$ and  $\alpha\le (1/17)^3$,
we obtain the following:
\begin{equation*}\label{U'size}
   \delta(G'')\ge 7m'/8+1.
\end{equation*}

Let $a_{2i}'\in \Gamma(a_{1i}, U_2')$,
$a_{1i}'\in \Gamma(a_{2i}, U_1')$ such that $a_{1i}'a_{2i}'\in E(G)$;
and
$b_{2i}'\in \Gamma(b_{1i}, U_2')$,
$b_{1i}'\in \Gamma(b_{2i}, U_1')$ such that $b_{1i}'b_{2i}'\in E(G)$.
We have the claim below.

\begin{CLA}\label{spanning ladder}
The balanced bipartite graph
$G''$ contains three vertex-disjoint ladders $Q_1$, $Q_2$, and $ Q_3$
spanning on $V(G'')$
such that the first rung of $Q_1$ is $x_1'x_2'$ and the last rung of $Q_1$
is $a_{11}'a_{21}'$, the first rung of $Q_2$ is $b_{11}'b_{21}'$ and the last rung of $Q_2$
is $a_{12}'a_{22}'$, the first rung of $Q_3$ is $b_{12}'b_{22}'$ and
the last rung of $Q_3$ is $y_1'y_2'$.
\end{CLA}

\pf Since  $\delta(G'')\ge 7m'/8+1>m'/2+6$,
$G''$ has a perfect matching  $M$ containing the following edges: $x_1'x_2', a_{11}'a_{21}', b_{11}'b_{21}', a_{12}'a_{22}', b_{12}'b_{22}', y_1'y_2'$.
We identify $a_{11}'$ and $b_{11}'$, $a_{21}'$ and $b_{21}'$, $a_{12}'$ and $b_{12}'$, and
$a_{22}'$ and $b_{22}'$   as vertices called $c_{11}'$,   $c_{21}'$,
 $c_{12}'$,   and $c_{22}'$, respectively.
Denote $G^*=G^*(U_1^*,U_2^*)$ as the resulting graph and
let $c_{11}'c_{21}', c_{12}'c_{22}'\in E(G^*)$.
Denote $M':=M-\{a_{11}'a_{21}', b_{11}'b_{21}', a_{12}'a_{22}', b_{12}'b_{22}'\}\cup \{c_{11}'c_{21}', c_{12}'c_{22}'\}$.
 Define an auxiliary graph $H'$ on $M'$ as
follows. If $xy, uv\in M'-\{c_{11}'c_{21}', c_{12}'c_{22}'\}$ with $x,u\in U_1'$ then $xy\sim_{H'} uv$ if and only if $x\sim_{G'} v$
and $y\sim_{G'} u$. For any $pq\in M'-\{c_{11}'c_{21}', c_{12}'c_{22}'\}$ with $p\in U_2'$,
 $pq\sim_{H'} c_{11}'c_{21}'$\,(resp.
 $pq\sim_{H'} c_{12}'c_{22}'$) if and only if $p\sim_{G'} a_{11}', b_{11}'$
 and $q\sim_{G'} a_{21}', b_{21}'$\,(resp. $p\sim_{G'} a_{12}', b_{12}'$
 and $q\sim_{G'} a_{22}', b_{22}'$).
Notice that  there is a natural one-to-one correspondence
 between ladders with rungs in $M'$ and paths in $H'$.
 Since  $\delta_{G^*}(U_1^*,U_2^*), \delta_{G^*}(U_2^*,U_1^*)\ge 3m'/4+1$,
we get $\delta(H')\ge m'/2+1$.
Hence $H'$ has a hamiltonian path starting with $x_1'x_2'$, ending with $y_1'y_2'$,
and having $c_{11}'c_{21}'$ and $c_{12}'c_{22}'$  as two  internal vertices.
The path with the vertex $c_{11}'c_{21}'$ replaced by $a_{11}'a_{21}'$ and $b_{11}'b_{21}'$,  and
with the vertex $c_{12}'c_{22}'$ replaced by $a_{12}'a_{22}'$ and $b_{12}'b_{22}'$
is corresponding to the required ladders in $G''$.
\qed

If $T\in \{T_1,T_2\}$, then
$$
H=x_1x_2Q_1L_1Q_2L_2Q_3y_1y_2\cup T
$$
is a spanning Halin subgraph of $G$.  Suppose now that $T\in \{T_3,T_4,T_5\}$
and $z$ is the pendent vertex.
Then $z\in V_1'\cup V_2'-W$ by Claim~\ref{subgraph_K}.
Suppose, w.l.o.g., that $z\in V_2'-W$. Then
by (iii) of Claim~\ref{subgraph_K} and $\delta(V_2', V_1')\ge (1/2-2\alpha_1/3)n$ from~\eqref{degrees},
we have  that $deg_G(z, U_1')\ge deg_G(z, V_1'-V(L_1\cup L_2)-V(T))\ge (1-\alpha_1-10\alpha_2)m>m/2+1$.
So $z$ has a neighbor on each side of the ladder $Q_1L_1Q_2L_2Q_3$, which has $m$
vertices on each side, and each side has at most
$m/2+1$ vertices from each partition of $U_1'$ and $U_2'$.
Let $H'$ be obtained from $x_1x_2Q_1L_1Q_2L_2Q_3y_1y_2\cup T$
by suppressing the degree 2 vertex $z$.  Then $H'$ is a Halin graph
such that there exists one side of
$Q_1L_1Q_2L_2Q_3$ with each vertex on it as a degree 3 vertex on a underlying tree of $H'$.
Let $z'$ be a neighbor of $z$ such that $z'$ has degree 3 in the
underlying tree of $H'$. Then
$$
H=x_1x_2Q_1L_1Q_2L_2Q_3y_1y_2\cup T\cup\{zz'\},
$$
is a spanning Halin subgraph of $G$.

\subsection{Proof of Theorem~\ref{non-extremal}}

In this section, we prove Theorem~\ref{non-extremal}.  In the first subsection,
we introduce the  Regularity Lemma,   the Blow-up Lemma,  and some related
results. Then we show that $G$ contains
a subgraph $T$ isomorphic to $T_1$ if $n$ is even and to $T_2$
if $n$ is odd. By showing that $G-V(T)$
contains a spanning ladder $L$ with its
first rung adjacent to the head link of $T$
and its last rung adjacent to the tail link of $T$,
we get a spanning Halin subgraph $H$ of $G$
formed by $L\cup T$.

\subsubsection{The Regularity Lemma and the Blow-up Lemma}

For any two disjoint non-empty vertex-sets $A$ and $B$ of a graph $G$, the \emph{density}
of $A$ and $B$ is the ratio $d(A,B):=\frac{e(A,B)}{|A|\cdot|B|}$. Let $\varepsilon$ and
$\delta$ be two positive real numbers. The pair $(A,B)$ is called $\ve$-regular if for
every $X\subseteq A$ and $Y\subseteq B$ with $|X|>\ve|A|$ and $|Y|>\ve|B|$,
$|d(X,Y)-d(A,B)|<\ve$ holds. In addition, if $\delta(A,B)>\delta |B|$ and
$\delta(B,A)>\delta|A|$, we say $(A,B)$ an $(\ve,\delta)$-super regular pair.

\begin{LEM}[\textbf{Regularity lemma-Degree form~\cite{Szemeredi-regular-partitions}}]\label{regularity-lemma}
For every $\ve>0$ there is an $M=M(\ve)$ such that if $G$ is any graph with $n$ vertices and $d\in[0,1]$ is
any real number, then there is a partition of the vertex-set $V(G)$ into $l+1$ clusters $V_0,V_1,\cdots,V_l$,
and there is a spanning subgraph $G'\subseteq G$ with the following properties.
\vspace{-2mm}
\begin{itemize}
  \item $l\le M$;
  \item $|V_0|\le \ve n$, all clusters $|V_i|=|V_j|\le \lceil\ve n\rceil$ for all $1\le i\ne j\le l$;
  \item $deg_{G'}(v)>deg_G(v)-(d+\ve)n$ for all $v\in V(G)$;
  \item $e(G'[V_i])=0$ for all $i\ge 1$;
  \item in $G'$, all pairs $(V_i,V_j)$ ($1\le i\ne j\le l$) are $\ve$-regular, each with a density either $0$ or greater than $d$.
 \end{itemize}
\end{LEM}

\begin{LEM}[\textbf{Blow-up lemma~\cite{Blow-up}}]\label{blow-up}
For every $\delta, \Delta, c>0$, there exists an $\ve=\ve(\delta, \Delta, c)$
and $\gamma=\gamma(\delta, \Delta, c)>0$ such that the following holds.
Let $(X,Y)$ be an $(\ve, \delta)$-super-regular pair with $|X|=|Y|=N$.
If a bipartite graph $H$ with $\Delta(H)\le \Delta$ can be embedded in
$K_{N,N}$ by a function $\phi$, then $H$ can be embedded in $(X,Y)$.
Moreover, in each $\phi^{-1}(X)$ and $\phi^{-1}(Y)$\,(the inverse image of
$X$ and $Y$, respectively), fix at most $\gamma N$
special vertices $z$, each of which is equipped  with a subset $S_z$
of $X$ or $Y$ of size at least $cN$.
The embedding of $H$ into $(X,Y)$ exists even if we restrict the image of
$z$ to be $S_z$ for all special vertices $z$.
\end{LEM}



Besides the above two lemmas, we also need the two lemmas below regarding regular pairs.

\begin{LEM}\label{regular-pair-large-degree}
If $(A,B)$ is an $\ve$-regular pair with density $d$, then for any $A'\subseteq A$ with
$|A'|\ge \ve |A|$, there are at most $\ve |B|$ vertices $b\in B$ such that
$deg(b, A')<(d-\ve)|A'|$.
\end{LEM}

\begin{LEM}[\textbf{Slicing lemma}]\label{slicing lemma}
Let $(A,B)$ be an $\ve$-regular pair with density $d$, and for some $\nu >\ve$, let $A'\subseteq A$
and $B'\subseteq B$ with $|A'|\ge \nu|A|$, $|B'|\ge \nu|B|$. Then $(A',B')$ is an $\ve'$-regular
pair of density $d'$, where $\ve'=\max\{\ve/\nu, 2\ve\}$ and $d'>d-\ve$.
\end{LEM}

%
%
%
%

\subsubsection{Finding subgraph $T$}

\begin{CLA}\label{find_K}
Let $n$ be a sufficient large integer and $G$  an $n$-vertex graph with $\delta(G)\ge (n+1)/2$.
If $\Delta(G)\le n-2$, then
$G$ contains a subgraph $T$ isomorphic to $T_1$
if $n$ is even and to $T_2$ if $n$ is odd.
\end{CLA}

\pf (The proof  can be much easier if uses the assumption that $G$ is not in Extremal Case 2, but we
show it here just using the conditions on minimum and maximum degrees.)
Suppose first that $n$ is even.
As $\Delta(G)\le n-2$, $G$ has
two vertices $x_1$ and $y$
such that $x_1y\not\in E(G)$.
Since $\delta(G)\ge (n+1)/2$, there exists $x\in \Gamma(x_1,y)$.
If $G$ has no independent set of size $n/2-1$,
we can find $x_2\in \Gamma(x_1,x)-\{y\}$
and $y_1,y_2\in \Gamma(y)-\{x,x_1,x_2\}$
such that $y_1y_2\in E(G)$.
Hence $G[\{x,y,x_1,x_2,y_1,y_2\}]$ contains a subgraph $T$ isomorphic to $T_1$.
So we assume $G$ has an independent set $S$ of size $n/2-1$.
Then $\delta(G-S)\ge 2$ and $\delta(S, V(G)-S)=|V(G)|-|S|$.
Let $x_1y_1$ and $yy_2$ be two independent
edges in $G-S$, and $x, y_1$ be any two
distinct vertices in $S$. Then $x,y_1\in \Gamma(x_1,y_1, y, y_2)$ and
$G[\{x,y,x_1,x_2,y_1,y_2\}]$ contains a subgraph $T$ isomorphic to $T_1$.



Then assume that $n$ is odd. Assume
first that $G$ has no independent set
of size $(n+1)/2-4$. We show in the first step that
$G$ contains a subgraph isomorphic to $K_4^-$\,($K_4$ with one edge removed).
Let $yz\in E(G)$. As $\delta(G)\ge (n+1)/2$, there exists
$y_1\in \Gamma(y,z)$. If there exists $y_2\in \Gamma(y,z)-\{y_1\}$, we are done.
Otherwise, $(\Gamma(y)-\{y_1,z\})\cap(\Gamma(z)-\{y_1,y\}) =\emptyset$.
As $\delta(G)\ge (n+1)/2$, $y_1$ is adjacent to
a vertex $y_2\in \Gamma(y)\cup \Gamma(z)-\{y_1,y,z\}$.
Assume $y_2\in \Gamma(z)-\{y_1,y\}$. Then $G[\{y,y_1,z,y_2\}]$
contains a copy of $K_4^-$. Choose $x\in \Gamma(y)-\{z,y_1, y_2\}$ and choose
an edge $x_1x_2\in G[\Gamma(x)-\{y,y_1, y_2,z\}]$.
Then $G[\{y,y_1,z,y_2, x, x_1, x_2\}]$ contains a subgraph $T$
isomorphic to $T_2$.  Let $S$ be a maximum  independent
set of $G$. So we have $|S|\ge (n+1)/2-4$.
Let $x,y_1,z\in S$. Then there exist
an edge $x_1x_2$ such that $x_1,x_2\in \Gamma(x,z)$,
and $y\in \Gamma(x,y_1)-\{x_1,x_2,z\}$.
If $|S|<(n+1)/2-2$, we
can find $y_2\in \Gamma(y,y_1)-\{x,x_1,x_2,z\}$.
Again, $G[\{y,y_1,z,y_2, x, x_1, x_2\}]$ contains a subgraph $T$
isomorphic to $T_2$. So $|S|\ge (n+1)/2-2$.
Assume $G-S$ does not contain any independent
set of size at least $(n+1)/2-6$.
For any $u,v,w\in S$, $\Gamma(u,v,w , V(G)-S)\ge (n+1)/2-4$
by Lemma~\ref{common-vertex}.
Let $x,y_1,z\in S$. As  $G-S$ does not contain any independent
set of size at least $(n+1)/2-6$,
there exist independent edges
$x_1x_2, yy_2$ such that
$x_1,x_2, y, y_2\in \Gamma(x,z,y_1, V(G)-S)$.
Then $G[\{y,y_1,z,y_2, x, x_1, x_2\}]$ contains a subgraph $T$
isomorphic to $T_2$. So
we assume that
$G-S$  contains an independent
set  of size at least $(n+1)/2-6$.
We take $Q$ to be a maximum one.
As $\delta(G-S)\ge 2$ and $Q$ is independent,
there exist two vertices  $x_1,y\in V(G)-S-Q$
such that both of them have at least
$(n+1)/13$ neighbors in $Q$.
Note that $\delta(S,V(G)-S)\ge (n+1)/2\ge |V(G)-S|-1$
and $S$ is a maximum independent set in $G$.
So there exist $x,z\in \Gamma(x_1, S)$
and $y_1\in \Gamma(y, S)-\{x,z\}$.
As both $x_1$ and $y$
have at least
$(n+1)/13$ neighbors in $Q$
and $\delta(S,V(G)-S)\ge (n+1)/2\ge |V(G)-S|-1$,
there exist $x_2\in \Gamma(x,x_1,z,Q)$
and $y_2\in \Gamma(y,y_1,Q)-\{x_2\}$.
So $G[\{y,y_1,z,y_2, x, x_1, x_2\}]$ contains a subgraph $T$
isomorphic to $T_2$.
%
\qed

Let $T$ be a subgraph of $G$ as given by Claim~\ref{find_K}.
Suppose the head link of $T$ is $x_1x_2$ and the tail link of $T$
is $y_1y_2$. Let $G'=G-V(T)$.  We show in next section
that $G'$ contains a spanning ladder with its first rung being adjacent to
$x_1x_2$ and its last rung being adjacent to $y_1y_2$.
 Let $n'=|V(G')|$. Then we have $\delta(G')\ge (n+1)/2-7\ge n'/2-4\ge
(1/2-\beta)n'$, where $\beta$ is the parameter  defined in the two extremal cases.

\subsubsection{Finding a spanning ladder of $G'$ with prescribed end rungs}


\begin{THM}\label{ladder}
Let $n'$ be a sufficiently large even integer and $G'$  an $n'$-vertex  subgraph of $G$
obtained by removing vertices in $T$, where $T\in \{T_1,T_2\}$ has head link $x_1x_2$
and tail link $y_1y_2$.
Suppose that $\delta(G')\ge (1/2-\beta)n'$ and $G=G[V(G')\cup V(T)]$
is in Non-extremal Case, then $G'$ contains a spanning ladder with its first rung adjacent to
$x_1x_2$ and its last rung adjacent to $y_1y_2$.
\end{THM}

\pf We fix the following sequence of parameters
\[
0<\ve\ll d\ll \beta \ll 1
\]
and specify their dependence as the proof proceeds.

Let $\beta$ be the parameter  defined in the two extremal cases.
Then we choose $d\ll \beta$ and  choose
\[
\ve=\frac{1}{4}\epsilon(d/2,3,d/4)
\]
following the definition of $\epsilon$ in the Blow-up Lemma.

Applying the Regularity Lemma to $G'$ with parameters $\ve$
and $d$, we obtain a partition of $V(G')$ into $l+1$ clusters
$V_0,V_1,\cdots, V_{l}$ for some $l \le M\le M(\ve)$, and a spanning subgraph
$G''$ of $G'$ with all described properties in the Regularity Lemma. In particular,
for all $v\in V(G')$,
\begin{equation}\label{G prime delta}
   deg_{G''}(v)>deg_{G'}(v)-(d+\ve)n'\ge (1/2-\beta-\ve-d)n'\ge (1/2-2\beta)n'
\end{equation}
provided that $\ve+d\le \beta$. On the other hand,
\begin{equation*}
    e(G'')\ge e(G')-\frac{(d+\ve)}{2}(n')^2>e(G')-d(n')^2
\end{equation*}
by $\ve<d$.

We further assume that $l=2k$ is even; otherwise, we eliminate the last cluster $V_{l}$
by removing all the vertices in this cluster to $V_0$. As a result, $|V_0|\le 2\ve n'$, and
\begin{eqnarray}\label{order_relation}
    (1-2\ve)n'\le l N=2kN\le n',
\end{eqnarray}
where $N=|V_i|$ for $1\le i\le l$.

For each pair $i$ and $j$ with $1\le i \ne j \le l$, we write $V_i\sim V_j$ if $d(V_i,V_j)\ge d$.
As in other applications of the Regularity Lemma, we consider the {\it reduced graph $G_r$},
whose vertex set is $\{1,2,\cdots, l\}$ and two vertices $i$ and $j$ are adjacent
if and only if $V_i\sim V_j$. From $\delta(G'')>(1/2-2\beta)n'$, we claim that $\delta(G_r)\ge (1/2-2\beta)l$.
Suppose not, and let $i_0\in V(G_r)$ be a vertex with $deg_{Gr}(i_0)<(1/2-2\beta)l$. Let $V_{i_0}$
be the cluster in $G$ corresponding to $i_0$. Then we have
\begin{equation*}
    (1/2-\beta)n'|V_{i_0}|\le |E_{G'}(V_{i_0},V-V_{i_0})|<(1/2-2\beta)l N|V_{i_0}|+2\ve n'|V_{i_0}| < (1/2-\beta)n'|V_{i_0}|.
\end{equation*}
This gives a contradiction by $l N\le n'$ from inequality~\eqref{order_relation}.

Let $A$ be a cluster of $G''$. We say $A$ is an {\it $(\ve, d)$-cluster\/} if
for any distinct cluster $B$ of $G''$ with $d(A,B)>0$, $(A,B)$ is an
$\ve$-regular pair with density at least $d$.
Let $x\in V(G')$ be a vertex and $A$ an $(\ve, d)$-cluster. We say $x$ is {\it typical} to $A$
if $deg(x,A)\ge (d-\ve)|A|$, and in this case, we write $x\sim A$.

\begin{CLA}\label{x1-y2}
Each vertex in $\{x_1,x_2,y_1,y_2\}$ is typical to at least $(1/2-2\beta)l$ clusters in $\{V_1, \cdots, V_l\}$.
\end{CLA}

\pf  Suppose on the contrary that there exists $x\in \{x_1,x_2,y_2,y_2\}$ such that $x$ is typical to
less than $(1/2-2\beta)l$ clusters in $\{V_1, \cdots, V_l\}$. Then we have $deg_{G'}(x)<(1/2-2\beta)l N+(d+\ve) n'\le (1/2-\beta)n'$
by $lN\le n'$ and $d+\ve \le \beta$.
\qed

Let $x\in V(G')$ be a vertex. Denote by  $\mathcal{V}_x$  the set of clusters to which $x$ typical.

\begin{CLA}\label{x1-x2}
There exist $V_{x_1}\in \mathcal{V}_{x_1}$ and $V_{x_2}\in \mathcal{V}_{x_2}$ such that
$d(V_{x_1}, V_{x_2})\ge d$.
\end{CLA}

\pf We show the claim by considering two cases based on the size of $| \mathcal{V}_{x_1}\cap  \mathcal{V}_{x_2}|$.

Case 1.  $| \mathcal{V}_{x_1}\cap  \mathcal{V}_{x_2}|\le 2\beta l$.

Then we have $| \mathcal{V}_{x_1}- \mathcal{V}_{x_2}|\ge (1/2-4\beta)l$ and $| \mathcal{V}_{x_2}- \mathcal{V}_{x_1}|\ge (1/2-4\beta)l$. We conclude that there is an edge between $\mathcal{V}_{x_1}- \mathcal{V}_{x_2}$
and $\mathcal{V}_{x_2}- \mathcal{V}_{x_1}$ in $G_r$.
For otherwise, let $\mathcal{U}$ be the union of clusters in
$\mathcal{V}_{x_1}\cap  \mathcal{V}_{x_2}$, $W=V_0\cup \mathcal{U}\cup V(T)$.
Let $V_1$ be the set of vertices contained in clusters in $ \mathcal{V}_{x_1}- \mathcal{V}_{x_2}$,
and $V_2$ be the set of vertices contained in clusters in $ \mathcal{V}_{x_1}- \mathcal{V}_{x_2}$.
Then $V_1$ and $V_2$ is a partition of $V(G)-W$.
Furthermore, $|W|\le 5\beta n$, $e(V_1, V_2)\le (d+\ve)n'|V_1|\le (d+\ve)n'(1+4\beta)lN\le \beta n^2$,
and $\delta(G[V_i])\ge \delta(G)-7-|W|-(d+\ve)n'\ge \delta(G)-|W|-\beta n$.
These imply that $W$ is an approximate vertex-cut of parameter $\beta$
with size at most $5\beta n$,
implying that $G$ is in Extremal Case 1.

Case 2.  $| \mathcal{V}_{x_1}\cap  \mathcal{V}_{x_2}|>2\beta l$.

We may assume that $\mathcal{V}_{x_1}\cap  \mathcal{V}_{x_2}$ is an independent set in $G_r$. For otherwise, we are done by finding an edge within
$\mathcal{V}_{x_1}\cap  \mathcal{V}_{x_2}$.
Also we
may assume that $E_{G_r}(\mathcal{V}_{x_1}\cap  \mathcal{V}_{x_2}, \mathcal{V}_{x_1}-\mathcal{V}_{x_2})=\emptyset$
and $E_{G_r}(\mathcal{V}_{x_1}\cap  \mathcal{V}_{x_2}, \mathcal{V}_{x_2}-\mathcal{V}_{x_1})=\emptyset$.
Since $\delta(G_r)\ge (1/2-2\beta)l$ and
$\delta_{G_r}(\mathcal{V}_{x_1}\cap  \mathcal{V}_{x_2}, \mathcal{V}_{x_1}\cup \mathcal{V}_{x_2})=0$,
we know that $l-|\mathcal{V}_{x_1}\cup \mathcal{V}_{x_2}|\ge (1/2-2\beta)l$.
Hence, $|\mathcal{V}_{x_1}\cup \mathcal{V}_{x_2}|=|\mathcal{V}_{x_1}|+|\mathcal{V}_{x_2}|-
|\mathcal{V}_{x_1}\cap  \mathcal{V}_{x_2}|\le (1/2+2\beta)l$.  This gives that $|\mathcal{V}_{x_1}\cap  \mathcal{V}_{x_2}|\ge |\mathcal{V}_{x_1}|+|\mathcal{V}_{x_2}|-(1/2+2\beta)l \ge (1/2-2\beta)l+(1/2-2\beta )l-(1/2+2\beta )l
\ge (1/2-6\beta )l$. Let $\mathcal{U}$ be the union of clusters in
$\mathcal{V}_{x_1}\cap  \mathcal{V}_{x_2}$.
Then $|\mathcal{U}|\ge (1/2-7\beta)n$
and $\Delta(G[\mathcal{U}])\le (d+\ve)n'\le \beta n$. This shows that $G$
is in Extremal Case 2.
\qed

Similarly, we have the following claim:
\begin{CLA}\label{y1-y2}
There exist $V_{y_1}\in \mathcal{V}_{y_1}-\{V_{x_1}, V_{x_2}\}$ and $V_{y_2} \in \mathcal{V}_{y_2}-\{V_{x_1}, V_{x_2}\}$ such that
$d(V_{y_1}, V_{y_2})\ge d$.
\end{CLA}

\begin{CLA}\label{hamiltonian_path}
The reduced graph $G_r$ has a hamiltonian path $X_1Y_1\cdots X_kY_k$ such that
$\{X_1, Y_1\}=\{V_{x_1}, V_{x_2}\}$ and $\{X_k, Y_k\}=\{V_{y_1}, V_{y_2}\}$.
\end{CLA}

\pf We contract the edges $V_{x_1}V_{x_2}$ and  $V_{y_1}V_{y_2}$ in $G_r$. Denote
the two new vertices as $V_x'$ and $V_y'$ respectively, and denote the resulting graph as $G_r'$.
Then we show that $G_r'$ contains a hamiltonian $(V_x', V_y')$-path. This  path is corresponding to
a required hamiltonian path in $G_r$.

To show $G_r'$ has a hamiltonian $(V_x', V_y')$-path, we need the following generalized version of a result
due to Nash-Williams~\cite{MR0284366} : Let $Q$ be a 2-connected graph of order $m$. If $\delta(Q)\ge \max\{(m+2)/3+1, \alpha(Q)+1\}$, then $Q$ is hamiltonian connected, where $\alpha(Q)$ is the size of a largest independent set of $Q$.

We claim  that $G_r'$ is $2\beta l$-connected.  Otherwise, let $S$ be a vertex-cut of $G_r'$ with $|S|<2\beta l$
and $\mathcal{S}$ the vertex set corresponding to $S$ in $G$. Since $\delta(G_r')\ge (1/2-2\beta)l-2$ and $|S|<2\beta l$,
we know that $G_r'-S$ has exactly two components.
Let $W=\mathcal{S}\cup V_0\cup V(T)$, $V_1$ the set of vertices contained in clusters
corresponding to vertices in one component of $G_r'-S$, and $V_2=V(G)-V_1-W$.
Then it is easy to check that $e(V_1,V_2)\le \beta n^2$ and $\delta(G[V_i])\ge \delta(G)-|W|-\beta n$.
Hence $W$ is an approximate vertex-cut of parameter $\beta$
with size at most $5\beta n$, showing that $G$ is in Extremal Case 1.
Since $n'=Nl+|V_0|\le (l+2)\ve n'$, we have that $l\ge 1/\ve -2\ge 1/\beta$. Hence, $G_r'$
is 2-connected.  As $G$ is not in Extremal Case 2, $\alpha(G_r')\le (1/2-7\beta)l$. By $\delta(G_r)\ge (1/2-2\beta)l$, we
have $\delta(G_r')\ge (1/2-2\beta)l-2\ge \max\{(l+2)/3+1, (1/2-7\beta)l+1\}$.  Thus, by the result on hamiltonian connectedness given above, we
know that $G_r'$ contains a hamiltonian $(V_x', V_y')$-path.
\qed


\begin{CLA}\label{super-regular}
For each $1\le i\le k$, there exist $X_i'\subseteq X_i$ and $Y_i'\subseteq Y_i$ such that
each of the following holds:
\begin{enumerate}[(1)]
  \item $|X_1'|\ge (1-\ve)|X_1|-1$, $|Y_k'|\ge (1-\ve)|Y_k|-1$, $|Y_1'|\ge (1-\ve)|Y_1|$,
$|X_k'|\ge (1-\ve)|X_k|$, and $|X_i'|\ge (1-\ve)|X_i|$, $2\le i\le k-1$;
  \item $(X_i', Y_i')$ is $(2\ve, d-3\ve)$-super-regular with density at least $d-\ve$;
  \item $|Y_1'|=|X_1'|+1$, $|X_k'|=|Y_k'|+1$, and $|X_i'|=|Y_i'|$, $2\le i\le k-1$; and
  \item for any  $A, B\in \{X_1', Y_1', \cdots, X_k', Y_k'\}$, if $d(A,B)>0$, then $(A,B)$
  is $2\ve$-regular with density at least $d-\ve$. Consequently, each $A$ is a $(2\ve, d-\ve)$ cluster.
\end{enumerate}
\end{CLA}

\pf
For each $1\le i\le k$, let
\begin{eqnarray*}
  X_i'' &=& \{x\in X_i\,|\, deg(x,Y_i)\ge (d-\ve)N\},\, \mbox{and} \\
  Y_i'' &=& \{y\in Y_i\,|\,deg(y,X_i)\ge (d-\ve)N\}.
\end{eqnarray*}
 If necessary,
we either take a subset $X_i'$ of $X_i''$
or take a subset $Y_i'$ of $Y_i''$
such that $|Y_1'|=|X_1'|+1$, $|X_k'|=|Y_k'|+1$, and $|X_i'|=|Y_i'|$ for $2\le i\le k-1$.
Since $(X_i,Y_i)$ is $\ve$-regular, we have
$|X_i''|, |Y_i''|\ge (1-\ve)N$.  This gives that $|X_1'|, |Y_k'|\ge (1-\ve)N-1$, $|Y_1'|\ge (1-\ve)N$,
$|X_k'|\ge (1-\ve)N$,
and $|X_i'|=|Y_i'|\ge (1-\ve)N$ for $2\le i\le k-1$.
As a result, we have $\deg(x, Y_i')\ge (d-2\ve)N$ for each $x\in X_i'$
and $deg(y, X_i')\ge (d-2\ve)N-1\ge (d-3\ve)N$ for each $y\in Y_i'$.
By the Slicing lemma\,(Lemma~\ref{slicing lemma}), $(X_i', Y_i')$ is $2\ve$-regular with density at least $d-\ve$.
Hence $(X_i', Y_i')$ is $(2\ve, d-3\ve)$-super-regular for each $1\le i\le k$.
The last assertion is again an application of the Slicing lemma.
\qed


For  $1\le i\le k$, we call each $X_i', Y_i'$ a {\it super-regularized cluster\,(sr-cluster)\/},
and call $X_i'$ and $Y_i'$ partners of each other and  write $P(X_i')=Y_i'$ and $P(Y_i')=X_i'$.
Denote $R=V_0\cup (\bigcup\limits_{i=1}^{k}((X_i\cup Y_i)-(X_i'\cup Y_i')))$.  Since $|(X_i\cup Y_i)-(X_i'\cup Y_i')|\le 2\ve N$
for $2\le i \le k-1$ and $|(X_1\cup Y_1)-(X_1'\cup Y_1')|, |(X_k\cup Y_k)-(X_k'\cup Y_k')|\le 2\ve N+1$,
we have $|R|\le 2\ve n+ 2k \ve N+2 \le 3\ve n'$.
As $n'$ is even and $|X_1'|+|Y_1'|+\cdots +|X_k'|+|Y_k'|$ is even,  we know $|R|$ is even.  We arbitrarily group  vertices in $R$ into
$|R|/2$ pairs. Given two vertices $u,v\in R$, we define a $(u,v)$-chain of length $2t$
as distinct sr-clusters $A_1, B_1, \cdots, A_t, B_t$ such that  $u\sim A_1\sim B_1\sim \cdots \sim A_t\sim B_t\sim v$
and each $A_j$ and $B_j$ are partners, in other words, $\{A_j,B_j\}=\{X_{i_j}', Y_{i_j}'\}$ for some $i_j\in \{1, \cdots, k\}$.
Recall here  $u\sim A_1$ means that $deg(u, A_1)\ge (d-3\ve)|A_1|$, and $A_1\sim B_1$ means that
the two vertices corresponding to $A_1$ and $B_1$ are adjacent in $G_r$.
We call such a chain of length $2t$ a {\it $2t$-chain\/}.

\begin{CLA}\label{absorbing-pre}
For each pair $(u,v)$ in $R$, we can find a $(u,v)$-chain of length at most 4 such that every sr-cluster is contained  in at most $d^2N/5$
chains.
\end{CLA}

\pf Suppose we have found chains for the first $m<2\ve n'$ pairs
of vertices in $R$ such that no sr-cluster  is contained in more than  $d^2N/5$ chains.
Let $\Omega$ be the set of all sr-clusters  that are contained  exactly in $d^2N/5$ chains. Then
\begin{eqnarray*}
  \frac{d^2N}{5}|\Omega| &\le & 4m <8\ve n' \le 8\ve\frac{2kN}{1-2\ve},
   \end{eqnarray*}
where the last inequality follows from (\ref{order_relation}).
Therefore,
\begin{eqnarray*}
  |\Omega| &\le &\frac{80k\ve}{d^2(1-2\ve)}\le \frac{80l \ve}{d^2}\le \beta l/2,
  \end{eqnarray*}
provided that $1-2\ve \ge 1/2$ and $80\ve \le d^2\beta /2$.

Consider now a pair  $(w,z)$ of vertices in $R$ which does not have a chain found so far, we want to find a
$(w,z)$-chain using sr-clusters not in $\Omega$.  Let $\mathcal{U}$
be the set of all  sr-clusters  to which $w$ typical  but not in $\Omega$, and
let $\mathcal{V}$ be the set of all  sr-clusters to which  $z$ typical but not in $\Omega$.
We claim that $|\mathcal{U}|,|\mathcal{V}|\ge (1/2-2\beta)l$.
To see this, we first observe that any vertex $x\in R$ is typical to
at least $(1/2-3\beta/2)l$  sr-clusters. For instead,
\begin{eqnarray*}
  (1/2-\beta)n' &\le &deg_{G'}(x) < (1/2-3\beta/2)lN+ (d-3\ve)lN+ 3\ve n', \\
   &\le& (1/2-3\beta/2+d)n' \\
   &<&  (1/2-\beta)n' \,\,(\mbox{provided that $d< \beta/2$ }),
\end{eqnarray*}
showing a contradiction. Since $|\Omega|\le \beta l/2$,
we have $|\mathcal{U}|,|\mathcal{V}|\ge (1/2-2\beta)l$.
Let $P(\mathcal{U})$ and $P(\mathcal{V})$ be the set of
the partners of clusters in $\mathcal{U}$
 and $\mathcal{V}$, respectively.
 By the definition of the chains, a cluster $A\in \Omega$
 if and only its partner $P(A)\in \Omega$.  Hence, $(P(\mathcal{U})\cup P(\mathcal{V}))\cap \Omega=\emptyset$.
 Notice also that each cluster has a unique partner, and so we have
 $|P(\mathcal{U})|=|\mathcal{U}|\ge (1/2-2\beta)l $ and $|P(\mathcal{V})|=|\mathcal{V}|\ge (1/2-2\beta)l $.


%
%

 If $E_{G_r}(P(\mathcal{U}), P(\mathcal{V}))\ne \emptyset$, then there exist two adjacent clusters
 $B_1\in P(\mathcal{U})$, $A_2\in P(\mathcal{V})$. If $B_1$ and $A_2$ are partners of each other,
 then $w\sim A_2\sim B_1\sim z$ gives a $(w,z)$-chain of length 2.
 Otherwise, assume $A_1=P(B_1)$ and $B_2=P(A_2)$, then $w\sim A_1\sim B_1\sim A_2\sim B_2\sim z$
 gives a $(w,z)$-chain of length 4.
 Hence we assume that $E_{G_r}(P(\mathcal{U}), P(\mathcal{V}))= \emptyset$.
 We may assume that $P(\mathcal{U})\cap P(\mathcal{V}) \ne \emptyset $.
 Otherwise, let $\mathcal{S}$ be the union of clusters
 contained in $V(G_r)-(P(\mathcal{U})\cup P(\mathcal{V}))$.
 Then $\mathcal{S}\cup R\cup V(T)$ with $|\mathcal{S}\cup R\cup V(T)|\le 4\beta n'+3\ve n'+7\le 5\beta n$
 \,(provided that $3\ve +7/n'<\beta$)
 is an approximate  vertex-cut of $G$, implying that $G$ is in Extremal Case 1.
As $E_{G_r}(P(\mathcal{U}), P(\mathcal{V}))= \emptyset$, any cluster in
$P(\mathcal{U})\cap P(\mathcal{V})$ is adjacent to at least $(1/2-2\beta)l$ clusters
in $V(G_r)-(P(\mathcal{U})\cup P(\mathcal{V}))$ by $\delta(G_r)\ge (1/2-2\beta)l$.
This implies that $|P(\mathcal{U})\cup P(\mathcal{V})|\le (1/2+2\beta)l$,
and thus $|P(\mathcal{U})\cap P(\mathcal{V})|\ge |P(\mathcal{U})|+|P(\mathcal{V})|-|P(\mathcal{U})\cup P(\mathcal{V})|\ge
(1/2-6\beta)l$. Then $P(\mathcal{U})\cap P(\mathcal{V})$
is corresponding to a subset $V_1$ of $V(G)$ such that $|V_1|\ge (1/2-6\beta)lN\ge (1/2-7\beta)n$
and $\Delta(G[V_1])\le (d+\ve)n'\le \beta n$. This implies that $G$ is in Extremal Case 2,
showing a contradiction.
\qed

By Claim~\ref{absorbing-pre}, each  vertex in $R$
is contained in a unique chain of length at most 4.
Let $Z$ be an sr-cluster, and $u\in R$
be a vertex.
We say  $u$ and $Z$ are {\it chain-adjacent\/} to each other if in the chain which  contains $u$,
$Z$ appears next to $u$.
For each sr-cluster  $Z\in \{X_1', Y_1', \cdots, X_k', Y_k'\}$,
let $R(Z)$ denote the set of vertices in $R$ such that each of
the vertices is chain-adjacent to $Z$.
Let $R_4(Z)=\{u\in R(Z)\,|\, u\,\mbox{is contained in a 4-chain}\}$,
and let $S_4(Z)$ denote the set of sr-clusters distinct from $Z$ such that
each of them is adjacent to the partner $P(Z)$ of $Z$ in a 4-chain which contains
$Z$. That is, for each $A\in S_4(Z)$, there exists $u\in R_4(Z)$ and $v\in R-R_4(Z)$
such that $u\sim Z\sim P(Z)\sim A\sim P(A)\sim v$ is a 4-chain.
If $Z\in \{X_1', \cdots, X_k'\}$,
then for each sr-cluster $A\in S_4(Z)$, let $c(A)$ denote the number of 4-chains
which contains $Z\sim P(Z)\sim A\sim P(A)$ as a sequence.
For each $A\in S_4(Z)$, choose $c(A)$ vertices in $A$
such that each of them has at least $(d-3\ve)|Z|>3d^2N/5$ neighbors in $P(Z)$. (Since $(P(Z), A)$ is $2\ve$-regular with density at least
$d-\ve$, we know that there are at least $(1-2\ve)|A|$ vertices in $A$ with this property by Lemma~\ref{regular-pair-large-degree}.)
Let $R'(Z)$ be the union of $R(Z)$ and the set of vertices chosen from $A\in S_4(P(Z))$ above, and let
$$
\omega(A)=\sum\limits_{A\in S_4(Z),\,\, Z\in\{X_1', \cdots, X_k'\}} c(A).
$$
Note that by the definitions,
$R'(Z)$ is only defined for sr-clusters $Z\in \{Y_1', \cdots, Y_k'\}$,
and $\omega(A)$ is defined only for sr-clusters $A\in \{X_1',\cdots, X_k'\}$.

%

%
\begin{CLA}\label{small-ladders1}
For each $i=1,2,\cdots, k$, each of the following holds.
\begin{enumerate}[(a)]
\item $|R(X_i')|\le d^2N/5$ and $|R'(Y_i')|\le d^2N/5$.
  \item $|R(X_i')-R_4(X_i')|=|R(Y_i')-R_4(Y_i')|$.
  \item $\omega(X_i')=|R_4(Y_{i}')|$.
  \item $|R'(Y_i')-R(Y_i')|=|R_4(X_i')|$.
\end{enumerate}
\end{CLA}

\pf By Claim~\ref{absorbing-pre}, each sr-cluster is contained in at most $d^2N/5$
chains, and a chain contains $X_i'$ if and only if it also contains $Y_i'$ by its definition.
 Since both $|R(X_i')|$ and $|R'(Y_i')|$ are bounded above by 
the number of chains which contain them, we have that $|R(X_i')|\le d^2N/5$ and $|R'(Y_i')|\le d^2N/5$.
By the definition of 2-chains, a vertex in $R$ is chain-adjacent to an sr-cluster $A$
in a 2-chain if and only if there exists another vertex in $R$ which is
chain-adjacent to the partner $P(A)$ of $A$.
Thus $|R(X_i')-R_4(X_i')|=|R(Y_i')-R_4(Y_i')|$.
By the definition, if $X_i'\in S_4(Z)$ for some sr-cluster $Z$,
then $c(X_i')$
is the number of 4-chains which contains $Y_i'\sim X_i'\sim P(Z)\sim Z$ as a sequence.
All of  such 4-chains is just the set of 4-chains in which 
 $Y_i'$ is chain-adjacent to a vertex in $R$.
Since each vertex in $R$ is contained in a unique chain, we
then have that $\omega(X_i')=|R_4(Y_{i}')|$.
Since each vertex in $R'(Y_i')-R(Y_i')$ is corresponding to a 4-chain
 in which $X_i'$ is chain-adjacent to a vertex in $R$, we have that  $|R'(Y_i')-R(Y_i')|=|R_4(X_i')|$.
\qed

%

\begin{CLA}\label{small-ladders}
For each $i=1,2,\cdots, k$,  there  exist  vertex-disjoint
ladders $L_{x}^i$, $L_y^i$
such that
\begin{enumerate}[(a)]
\item $R(X_i')\subseteq V(L_{x}^{i})\subseteq  R(X_i')\cup  X_{i}'\cup Y_{i}'$
and $R'(Y_i')\subseteq V(L_{y}^{i})\subseteq  X_{i}'\cup Y_{i}'\cup R'(Y_i')$;
\item $|(V(L_x^i)\cup V(L_y^i))\cap X_i'|=4|R(X_i')|+3|R(Y_i')|+3|R_4(X_i')|-2$ and $|(V(L_x^i)\cup V(L_y^i))\cap Y_i'|=4|R(Y_i')|+4|R_4(X_i')|+3|R(X_i')|-2$; and
  \item the vertices on the first  and last rungs of each of  $L_{x}^i$ and  $L_y^i$ are contained in $X_i'\cup Y_i'$.
\end{enumerate}
%
\end{CLA}

\pf
Notice that by Claim~\ref{super-regular}, $(X_i', Y_i')$ is $2\ve$-regular with density at least $d-\ve$.
Let $R(X_i')=\{x_1, \cdots, x_r\}$.
For each $j$, $1\le j\le r$, since $|\Gamma(x_j,X_i')|\ge (d-3\ve)|X_i'|>2\ve|X_i'|$,
by Lemma~\ref{regular-pair-large-degree}, there exists a vertex set $B_j\subseteq Y_i'$
with $|B_j|\ge (1-2\ve)|Y_i'|$ such that for each $b_1\in B_j$, $deg(b_1,\Gamma(x_j,X_i'))\ge (d-3\ve)|\Gamma(x_j,X_i')|>4|R(X_i')|$.
If $r\ge 2$, for $j=1,\cdots, r-1$, by Lemma~\ref{regular-pair-large-degree}, there also exists a vertex set $B_{j,j+1}\subseteq Y_i'$
with
 $|B_{j,j+1}|\ge (1-4\ve)|Y_i'|$
such that
 for each $b_2\in B_{j,j+1}$, we have $deg(b_2,\Gamma(x_j,X_i'))\ge (d-3\ve)|\Gamma(x_{j},X_i')|>4|R(X_i')|$ and  $deg(b_2,\Gamma(x_{j+1},X_i'))\ge (d-3\ve)|\Gamma(x_{j+1},X_i')|>4|R(X_i')|$.
When $r\ge 2$,
since $|B_j|, |B_{j,j+1}|, |B_{j+1}|\ge (d-3\ve)|Y_i'|>2\ve |Y_i'|$,
there is a set $A\subseteq X_i'$ with $|A|\ge (1-6\ve)|X_i'|\ge |R(X_i')|$ such that for each $a\in A$,
  $deg(a,B_j)\ge (d-3\ve)|B_j|$, $deg(a,B_{j,j+1})\ge (d-3\ve)|B_{j,j+1}|$ and $deg(a,B_{j+1})\ge (d-3\ve)|B_{j+1}|$.
Notice that $(d-3\ve)|B_j|, (d-3\ve)|B_{j,j+1}|, (d-3\ve)|B_{j+1}|\ge (d-3\ve)(1-4\ve)|Y_i'|>3|R(X_i')|$.
Hence we can choose distinct  vertices $u_1, u_2, \cdots, u_{r-1}\in A$
such that $deg(u_j, B_{j}), deg(u_j, B_{j,j+1}), deg(u_j, B_{j+1})\ge 3|R(X_i')|$.
Then we can choose distinct vertices  $y_{23}^j\in \Gamma(u_j, B_{j}), z_{j}\in \Gamma(u_j, B_{j,j+1})$
and $y_{12}^{j+1}\in \Gamma(u_j, B_{j+1})$ for each $j$,
and choose distinct and unchosen vertices $y_{12}^1\in B_1$ and $y_{23}^r\in B_r$.
Finally, as
for each vertex $b_1\in B_j$, we have $deg(b_1,\Gamma(x_j,X_i'))>4|R(X_i')|$
and  for each vertex $b_2\in B_{j,j+1}$, we have $deg(b_2,\Gamma(x_j,X_i')), deg(b_2,\Gamma(x_{j+1},X_i'))>4|R(X_i')|$,
we can choose $x_{j1}, x_{j2}, x_{j3}\in \Gamma(x_j, X_i')-\{u_1,\cdots, u_{r-1}\}$
such that $y_{12}^j\in \Gamma(x_{j1}, x_{j2}, Y_i')$, $y_{23}^j\in \Gamma(x_{j2}, x_{j3}, Y_i')$,
and $z_j\in \Gamma(x_{i3}, x_{i+1,1}, Y_i')$.
(When $i\ge 2$, we choose all these vertices such that they are not used by existing ladders. The possibility of doing this 
is guaranteed by the degree conditions and the small sizes of the existing ladders.)
Let $L_{x}^i$ be the graph with
$$
V(L_{x}^i)=R(X_i')\cup \{x_{i1}, x_{i2}, x_{i3},y_{12}^i, y_{23}^i, z_i, u_i, x_{r1}, x_{r2}, x_{r3}, y_{12}^r, y_{23}^r \,|\, 1\le i\le r-1 \}\quad \mbox{and}
$$
$E(L_{x}^i)$ consisting of the edges $x_rx_{r1}, x_rx_{r2},x_rx_{r3}, y_{12}^rx_{r1}, y_{12}^rx_{r2},y_{23}^rx_{r2},y_{23}^rx_{r3}$
and the edges indicated below for each $1\le i\le r-1$:
$$
x_i\sim x_{i1},x_{i2}, x_{i3};\, y_{12}^i\sim x_{i1}, x_{i2};\, y_{23}^i\sim x_{i2},x_{i3}; \,  z_i\sim x_{i3}, x_{i+1,1}; \, u_i\sim x_{i3}, x_{i+1,1}, z_i.
$$
It is easy to check that $L_{x}^i$ is a ladder  spanning on $R(X_i')$, $4|R(X_i')|-1$
vertices from $X_i'$ and $3|R(X_i')|-1$ vertices from $Y_i'$.
Similarly, we
can find a ladder  $L_{y}^i$ spanning on $R'(Y_i')$, $4|R'(Y_i')|-1$
vertices from $Y_i'$ and $3|R'(Y_i')|-1$ vertices from $X_i'$. The constructions of ladders $L_x^i$ and $L_y^i$
verify both of statements (a) and (c). The statement (b) is seen
by the construction of the ladders and (d) of Claim~\ref{small-ladders1} which says that
 $|R'(Y_i')|=|R(Y_i')|+|R_4(X_i')|$.
\qed

For each $i=1,2,\cdots, k-1$, let $X_i^{**}=X_i'-V(\bigcup_{i=1}^k(L_{x}^i\cup L_{y}^i))$ and
$Y_i^{**}=Y_i'-V(\bigcup_{i=1}^k(L_{x}^i\cup L_{y}^i))$. Using Lemma~\ref{regular-pair-large-degree}, for $i\in \{1,\cdots, k-1\}$,
choose $y_i^*\in Y_i^{**}$ such that $|A_{i+1}|\ge dN/4$,
where $A_{i+1}:=X_{i+1}^{**}\cap \Gamma(y_i^*)$.
This is possible, as $(Y_i^{**}, X_{i+1}^{**})$ is $4\ve$-regular
with density at least $d-3\ve$.
\,(Applying Slicing lemma based on $(Y_i', X_{i+1}')$.) Similarly,
choose $x_{i+1}^*\in A_{i+1}$ such that $|D_{i}|\ge dN/4$,
where $D_{i}:=Y_i^{**}\cap \Gamma(x_{i+1}^*)$.
Let $S=\{y_{i}^*, x_{i+1}^*\,|\, 1\le  i \le k-1\}$, and
let  $X_i^*=X_i^{**}-S$  and $Y_i^*=Y_i^{**}-S$. We have the following holds.

\begin{CLA}\label{final-super-pair}
For each $i=1,2,\cdots, k$, $|X_i^*|=|Y_i^*|$ and $(X_i^*, Y_i^*)$ is $(4\ve, d/2)$-super-regular.

\end{CLA}

\pf
We show that   $|X_i^*|=|Y_i^*|$ for each $i$, $1\le i\le k$.
Since $|Y_1'|=|X_1'|+1$, $|X_k'|=|Y_k'|+1$, and $|X_i'|=|Y_i'|$ for $2\le i\le k-1$, and
$|X_1^{**}|=|X_1^*|$, $|Y_k^{**}|=|Y_k^*|$, and $|X_i^{**}|=|X_i^*|-1$, $|Y_j^{**}|=|Y_j^*|-1$
for $2\le i\le k, 1\le j\le k-1$, it suffices to show that
$|X_i'\cap V(\bigcup_{i=1}^k(L_{x}^i\cup L_{y}^i))|=|Y_i'\cap V(\bigcup_{i=1}^k(L_{x}^i\cup L_{y}^i))|$.
This is clear by (b) of  Claim~\ref{small-ladders} and Claim~\ref{small-ladders1}.
As
\begin{eqnarray*}
\begin{array}{lll}
   & |X_i'\cap V(\bigcup_{i=1}^k(L_{x}^i\cup L_{y}^i))|= 4|R(X_i')|+3|R(Y_i')|+3|R_4(X_i')|-2+\omega(X_i') &  \\
  = & 4|R(X_i')-R_4(X_i')|+3|R(Y_i')-R_4(Y_i')|+7|R_4(X_i')|+3|R_4(Y_i')|-2+\omega(X_i') &  \\
  = & 7|R(X_i')-R_4(X_i')|+7|R_4(X_i')|+4|R_4(Y_i')|-2, &
\end{array}
 \end{eqnarray*}
and
\begin{eqnarray*}
\begin{array}{lll}
  &|Y_i'\cap V(\bigcup_{i=1}^k(L_{x}^i\cup L_{y}^i))|= 3|R(X_i')|+4|R(Y_i')|+4|R_4(X_i')|-2 &\\
   =& 3|R(X_i')-R_4(X_i')|+4|R(Y_i')-R_4(Y_i')|+7|R_4(X_i')|+4|R_4(Y_i')|-2&\\
   =& 7|R(X_i')-R_4(X_i')|+7|R_4(X_i')|+4|R_4(Y_i')|-2.&
    \end{array}
 \end{eqnarray*}

Since $|R(X_i')|, |R'(Y_{i}')|\le d^2N/5$ for each $i$,  by the first part of argument, 
$|X_i'\cap V(\bigcup_{i=1}^k(L_{x}^i\cup L_{y}^i))|\le 4|R(X_i')|+4|R'(Y_i')|-2\le 2d^2N-2$ and 
$|Y_i'\cap V(\bigcup_{i=1}^k(L_{x}^i\cup L_{y}^i))|\le 4|R(X_i')|+4|R'(Y_i')|-2\le 2d^2N-2$. 
Thus $|X_i^*|, |Y_i^*|\ge (1-\ve-2d^2)N$.  
As $\ve, d \ll 1$, we
can assume that $1-\ve-2d^2<1/2$. Thus, by Slicing lemma based on the $2\ve$-regular
pair $(X_i', Y_i')$, we know that
$(X_i^*, Y_i^*)$ is $4\ve$-regular. 
Recall from Claim~\ref{super-regular} that $(X_i', Y_i')$ is $(2\ve, d-3\ve)$-super-regular,
we know that
for each $x\in X_i^*$, $deg(x, Y_i^*)\ge (d-3\ve-2d^2)|Y_i^*|>
d|Y_i^*|/2$. Similarly, we have
for each $y\in Y_i^*$, $deg(y, X_i^*)\ge
d|X_i^*|/2$. Thus $(X_i^*, Y_i^*)$ is $(4\ve, d/2)$-super-regular.
\qed


For each $i=1,2,\cdots, k-1$,
now set $B_{i+1}:=Y_{i+1}^*\cap \Gamma(x_{i+1}^*)$
and $C_i:=X_i^*\cap \Gamma(y_i^*)$.
Since $(X_i^*, Y_i^*)$ is $(4\ve, d/2)$-super-regular,
we have $|B_i|, |C_i|\ge d|X_i^*|/2>d|X_i^*|/4$.
Recall from Claim~\ref{hamiltonian_path} that
$\{X_1, Y_1\}=\{V_{x_1}, V_{x_2}\}$ and $\{X_k, Y_k\}=\{V_{y_1}, V_{y_2}\}$. We assume, w.l.o.g., that $X_1=
V_{x_1}$ and $X_k=V_{y_1}$.
Let $A_{1}=X_1^*\cap \Gamma(x_1)$, $B_{1}=Y_1^*\cap \Gamma(x_2)$, $C_{k}=X_k^*\cap \Gamma(y_1)$, and $D_{k}=Y_k^*\cap \Gamma(y_2)$.  Since $deg(x_1, X_1)\ge (d-\ve)N$, we have
$deg(x_1, X^*_1)\ge (d-\ve-2\ve-2d^2)N\ge d|X_1^*|/4$,
and thus $|A_1|\ge d|X_1^*|/4$.
Similarly, we have $|B_1|, |C_k|, |D_k|\ge d|X_1^*|/4$.
For each $1\le i\le k$, we assume that $L_{x}^i=a_1^ib_1^i-L_{x}^i-c_1^id_1^i$ and
$L_{y}^i=a_2^ib_2^i-L_{y}^i-c_2^id_2^i$,
where
$a_j^i, c_j^i\in Y_i'\subseteq Y_i$ and $b_j^i, d_j^i\in X_i'\subseteq X_i$  for $j=1,2$.
For $j=1,2$, let $A_j^i=X_i^*\cap \Gamma(a_j^i)$,
$C_j^i=X_i^*\cap \Gamma(c_j^i)$,
$B_j^i=Y_i^*\cap \Gamma(b_j^i)$,  and
$D_j^i=Y_i^*\cap \Gamma(d_j^i)$.
Since $(X_i', Y_i')$ is $(2\ve, d-3\ve)$-super-regular,
for $j=1,2$, we have
$|\Gamma(a_j^i, X_i')|, |\Gamma(c_j^i, X_i')|\ge (d-3\ve)|X_i'|$
and
$|\Gamma(b_j^i, Y_i')|, |\Gamma(d_j^i, Y_i')|\ge (d-3\ve)|Y_i'|$.
Thus, we have
$|A_j^i|, |B_j^i|, |C_j^i|, |D_j^i| \ge (d-3\ve)|X_i'|-2d^2N\ge d|X_i^*|/4=d|Y_i^*|/4$.

We now apply the Blow-up lemma on $(X_i^*, Y_i^*)$ to find a spanning ladder $L^i$ with its first and last rungs being contained
in $A_i\times B_i$ and $C_i\times D_i$, respectively,
and for $j=1,2$, its $(2j)$-th and $(2j+1)$-th rungs
being contained in $A_j^i\times B_j^i$ and
$C_j^i\times D_j^i$, respectively. We can then insert
$L_{x}^i$ between the 2nd and 3rd rungs of $L^i$ and
$L_{y}^i$ between the 4th and 5th rungs of $L^i$
to obtained a ladder $\mathcal{L}^i$ spanning on $X_i\cup Y_i-S$.  Finally, $\mathcal{L}^1y_1^*x_2^*\mathcal{L}^2\cdots y_{k-1}^*x_k^* \mathcal{L}^k$ is a spanning ladder of $G'$
with its first rung adjacent to $x_1x_2$ and its last rung adjacent to $y_1y_2$.

The proof is now complete.
\qed

%


%
%
%
%
%
%
%
%
%
%
%
%

\bibliographystyle{plain}
\bibliography{SSL-BIB}

\end{document}